\documentclass[preprint,11pt]{article}
\usepackage{amssymb}
\usepackage{mathrsfs}
\usepackage{amsfonts}
\usepackage{graphicx}
\usepackage{colortbl,dcolumn}
\usepackage{tabularx}
\usepackage{amsmath}
\usepackage{psfrag}
\usepackage{booktabs}
\usepackage{array}
\usepackage{cite}
\numberwithin{equation}{section}
\usepackage[top=1.2in, bottom=1.2in, left=0.9in, right=0.9in]{geometry}

\newcommand{\R}{\mathbb{R}}
\newcommand{\N}{\mathbb{N}}

\newcommand{\dd}{\text{d}}
\newtheorem{theorem}{Theorem}[section]

\newtheorem{assumption}[theorem]{Assumption}

\newtheorem{lemma}[theorem]{Lemma}

\begin{document}

%
 \title {Strong convergence rates of a fully discrete scheme for the Cahn-Hilliard-Cook equation
 \footnote{R.Q. was supported  by NSF of China (No. 11701073).  M.C. and X.W. were supported by NSF of China
 (Nos. 11971488, 12071488) and NSF of Hunan province (No. 2020JJ2040).
%
}
}

\author{
Ruisheng Qi$\,^\text{a}$,
\quad Meng Cai$\,^\text{b}$,
\quad Xiaojie Wang$\,^\text{b}$,
\\
\footnotesize $\,^\text{a}$ School of Mathematics and Statistics, Yancheng Teachers University, Yancheng, China\\
\footnotesize qiruisheng123@126.com\\
\footnotesize $\,^\text{b}$ School of Mathematics and Statistics, HNP-LAMA, Central South University, Changsha, China\\
\footnotesize csumathcai@csu.edu.cn, \quad
\footnotesize x.j.wang7@csu.edu.cn \\
}
\maketitle
\begin{abstract}\hspace*{\fill}\\
  \normalsize


The first aim of this paper is to examine existence, uniqueness and regularity for the Cahn-Hilliard-Cook (CHC) equation
in space dimension $d\leq 3$. By applying a spectral Galerkin method to the infinite {\color{black}{dimensional}} equation,
we elaborate the well-posedness and regularity of the finite dimensional approximate problem.
The key idea lies in transforming the stochastic problem {\color{black}{with additive noise}} into an equivalent random equation.
The regularity of the solution to the equivalent random equation is obtained, in one dimension, with the aid of
the Gagliardo-Nirenberg inequality and done in two and three dimensions, by  the energy argument.
Further, the approximate solution is shown to be strongly convergent to the unique mild solution of the original CHC equation,
whose spatio-temporal regularity can be attained by similar arguments.
In addition, a fully discrete approximation of such problem is investigated, performed
by the spectral Galerkin method in space and the backward Euler method in
time. The previously obtained regularity results of the problem help us to identify
strong convergence rates of the fully discrete scheme.

\end{abstract}

\section{Introduction}
During the last decades, there have been a large number of works devoted to numerical approximations of
stochastic partial differential equations (SPDEs),
see monographs \cite{Jentzen2011Taylor,Kruse2014Strong,Lord2014An} and references therein.
{\color{black}
Pioneering works have focused on the strong and weak convergence of numerical methods for SPDEs 
with globally Lipschitz continuous nonlinearities.
In the last decade, some techniques were proposed in some works to handle numerical approximations of 
SPDEs with non-globally Lipschitz continuous nonlinearities 
\cite{becker2017strong,becker2019strong,brehier2018strong,brehier2019analysis,brehier2018weak,feng2017finite,
liu2020strong,liu2020strong-multiplicative,cui2020weak,cui2019strong,
gyongy2016convergence,hutzenthaler2018strong,jentzen2019strong,jentzen2018exponential,kovacs2018discretisation,
kovacs15backward,Majee2017optimal,campbell2018adaptive,wang2020efficient,cai2021weak,liu2003convergence,
kovacs2014erratum,hutzenthaler2014perturbation,feng2020fully,
furihata2018strong,kovacs2011finite,ruisheng2018error,cui2021strong,cui2020absolute}. 
But not all problems of interest are covered and there remains a lot of work to do.}
A typical SPDE model with non-globally Lipschitz coefficients is the Allen-Cahn type SPDE,
which has been recently numerically studied  by many authors, e.g.,
\cite{becker2017strong,becker2019strong,brehier2018strong,brehier2019analysis,brehier2018weak,feng2017finite,
liu2020strong,liu2020strong-multiplicative,cui2020weak,cui2019strong,
gyongy2016convergence,hutzenthaler2018strong,jentzen2019strong,jentzen2018exponential,kovacs2018discretisation,
kovacs15backward,Majee2017optimal,campbell2018adaptive,wang2020efficient,cai2021weak,liu2003convergence}.
%
As another prominent SPDE model with non-globally Lipschitz coefficients,  the stochastic Cahn--Hilliard equations
{\color{black}{are also numerically investigated by many authors very recently}}
(see \cite{kovacs2014erratum,hutzenthaler2014perturbation,feng2020fully,
furihata2018strong,kovacs2011finite,ruisheng2018error,cui2021strong,cui2020absolute,li2016unconditionally,
hong2022finite,hong2022convergence}).
The present paper aims to further carry out  theoretical and numerical analysis of such equation.

Let $D$ be a bounded open set of $\mathbb{R}^d$, $d=1,2,3$ with   smooth boundary and let $H=L^2(D;\mathbb{R})$ be the real separable Hilbert space endowed with the usual inner product  and  norm and 
{\color{black}{
$\dot{H}:=\big\{v\in H:\int_Dv\,\dd x=0\big\} 
$.}}
This article is concerned with the following stochastic Cahn--Hilliard equation perturbed by additive noise,
\begin{align}\label{eq:stochastic-allen-Hilliard-equation}
\begin{split}
\left\{
\begin{array}{ll}
\dd u - \Delta w\dd t= \dd W(t),&\; \text{ in }\; D\times (0,T],
\\
w=-\Delta u+f(u), &\; \text{ in } \; D\times (0,T],
\\
\frac{\partial u}{\partial n}=\frac{\partial w}{\partial n}=0,& \; \text{ in } \;\partial D\times(0,T],
\\
u(0,x)=u_0,&\; \text{ in } \; D,
\end{array}
\right.
\end{split}
\end{align}
where $f(s)=s^3-s, s\in \mathbb{R}$. {\color{black}{That}} SPDE driven by additive noise is also called Cahn--Hilliard--Cook (CHC) equations in literature.
Following the framework of \cite{da1996stochastic}, we can rewrite the above problem as an abstract equation described by
\begin{align}\label{eq:abstract-CHC-form}
\,\dd X(t)
+
A(AX(t)+F(X(t)))\,\dd t
=
\,\dd W(t),\; X(0)=X_0,
\end{align}
where $A:D(A)\subset \dot{H}\rightarrow \dot{H}$ is the Neumann Laplacian and $-A^2$ generates an analytic semigroup $E(t)$ on $\dot{H}$.
Similarly as in  \cite{furihata2018strong,kovacs2011finite}, $\{W (t)\}_{t\geq 0}$ is assumed to be a $Q$-Wiener process in $H$ with respect to  a filtered probability space $(\Omega,\mathcal{F},\mathbb{P},\{\mathcal{F}_t\}_{t\geq0})$.
The nonlinear mapping $F$ is supposed to be  a Nemytskij operator, given by $F(u)(x)=f(u(x))$, $x\in D$.

The deterministic version of such equation has been extensively studied as a well-known model of a spinodal decomposition for a
binary mixture \cite{Cahn1958Free}. It can be also used to describe the diffusive process of populations and an oil film spreading over a
solid surface \cite{Cho1997The,Edwards1987Mathematical}.
For the stochastic version, one can consult, e.g., \cite{cardon2001cahn,elezovic1991stochastic,da1996stochastic,cui2019global,
cui2020absolute,antonopoulou2016existence}
for the existence, uniqueness and regularity results.
The first goal of this paper is to provide further regularity results for the mild solution to \eqref{eq:abstract-CHC-form}.
{\color{black}{Let $P$ be an orthogonal operator from $H $ to $\dot{H}$.}}  
Under further assumptions specified later, particularly including
\begin{align}\label{eq:condition-Q-A}
\|{\color{black}{A^{\frac{\gamma-2}2}P Q^{\frac12}\|_{{\color{black}{\mathcal{L}}_2(H)}}
<\infty,
\;\text{ for some } \; \gamma \in \big(\tfrac d 2, 4 \big],}}
\end{align}
Theorem \ref{them:regulairty-mild-solution} shows that the underlying  equation \eqref{eq:abstract-CHC-form}
admits a unique mild solution $X(t)$, defined by \eqref{them:eq-mild-solution-stochastic-equation},
which enjoys the following spatio-temporal regularities
\begin{align}\label{eq:sptial-regulaity}
X\in L^\infty\big([0,T]; L^p(\Omega;{\color{black}{H^\gamma}})\big),\;\forall p\geq 1,
\end{align}
and for $0\leq s<t\leq T$,
\begin{align}\label{eq:temporal-regulaity}
\|X(t)-X(s)\|_{L^p(\Omega;{\color{black}{H^\beta}})}
\leq
C(t-s)^{\min\{\frac12,\frac{\gamma-\beta}4\}},\,\forall\beta\in[0,\gamma].
\end{align}

In the following we compare findings in this article with existing regularity results in the literature
and also illustrate how \eqref{eq:sptial-regulaity} and \eqref{eq:temporal-regulaity} can be established. 
Under  the assumption that $A$ and $Q$ commute and
$
{\color{black}{
\mathrm{Tr}
( A^{ \gamma - 1} Q ) < \infty
}}
$
for some
$
\gamma>0,
$
 it was shown in \cite{da1996stochastic} that the solution belongs to $C([0, T], H)$
almost surely. Further, under  the assumption that $A$ and $Q$ commute and
$
\|A^{\frac12}Q^{\frac12}\|_{\mathcal{L}_2}^2=\mathrm{Tr}(AQ) < \infty,
$
the authors in \cite{da1996stochastic} also show $X(t)\in C([0,T],\dot{H}^1)$ and $\mathbb{E}[J(X(t))]<\infty$,
where $J(\cdot)$ is a  Lyapunov functional, given by
$
J(u)=\frac12\|\nabla u\|^2
+
\int_D\Phi(u)\,\dd x,\; \Phi'(s) = F ( s ), s \in \R.
$
In \cite{kovacs2011finite}, the authors  prove similar results  but under the weaker condition that $
\mathrm{Tr}(AQ) < \infty
$.
The key idea of \cite{da1996stochastic,kovacs2011finite} was to apply the It\^{o} formula to  $J(\cdot)$, which requires that the noise process $W(t)$ satisfies the condition $\|A^{\frac{\gamma-2}2}Q^{\frac12}\|_{\mathcal{L}_2}<\infty$ for $\gamma\geq 3$.
Although the analysis in \cite{da1996stochastic,kovacs2011finite} works quite well for CHC equation with smoother noise cases, 
it cannot be generalized to the rougher noise cases,{\color{black}{ especially in dimensions two and three.  
{\color{black}{
We would like to mention the references \cite{cui2021strong,hong2022finite,antonopoulou2016existence,cui2019global} 
on the global well-posedness and regularity estimates of CHC equation driven by the space-time white noise with $d =1$.}}
We emphasize that our spatial-temporal regularity result \eqref{eq:sptial-regulaity}-\eqref{eq:temporal-regulaity} 
is new in  dimensions two and three.}}

We now illuminate our approach to analyze the regularity of \eqref{eq:abstract-CHC-form}
which overcomes the difficulty caused by the spatially rough noise. {\color{black}{It is worthwhile to point out that this is possible 
because the noise is additive and this is a standard approach in that context.}}
{\color{black}{By introducing $Y(t) = X(t) - Z(t)$ with $Z(t) : = E(t)(I-P)X_0+\int_0^tE(t-s)\,\dd W(s)$,}}
we can reformulate \eqref{eq:abstract-CHC-form} as
\begin{align}\label{eq:solution-two-part-mild-solution}
\dot{Y}(t)+A\big(AY(t)+F(Y(t)+Z(t))\big)=0,\;{\color{black}{Y(0)=PX_0}}.
\end{align}
Such a problem can be viewed as a random PDE.
{\color{black}{By the same argument of the proof of  \cite[Theorem 3.1]{larsson2011finite}, it is easy to show  that
under the assumption $\|A^{\frac{\gamma-2}2}PQ^{\frac12}\|_{\mathcal{L}_2(H)}
<\infty$,\;for $\gamma\geq 0$ and $p\geq 1$,
 $Z(t)$ enjoys the following spatial regularity,}}
 \begin{align}
{\color{black}{\|Z(t)\|_{L^p(\Omega;H^\gamma)}\leq C(\|A^{\frac{\gamma-2}2}PQ^{\frac12}\|_{\mathcal{L}_2(H)}
+
\|X_0\|_{L^p(\Omega;H^\gamma)}),\;t>0.}}
 \end{align}
Hence it suffices to analyze the spatio-temporal regularities of the solution $Y( \cdot )$ to \eqref{eq:solution-two-part-mild-solution},
which relies on  Galerkin's method and energy arguments in this paper. The difficulty lies on showing $\mathbb{E}\big[\sup_{s\in[0,T]}\|Y(s)\|_{L^6}^p\big]<\infty$. For dimension one, this {\color{black}{can be  shown based on the idea of \cite{cui2021strong}}} and for dimensions two and three, it can be done by multiplying both sides of \eqref{eq:solution-two-part-mild-solution} by $Y_t$ and some other manipulations. Hence our approach works well for CHC equation with rough noise.

The second aim of this article is devoted to the error estimate of the fully discrete approximation of \eqref{eq:abstract-CHC-form}. Let $k=T/M$, $M\in \mathbb{N}$ be a uniform time step-size and {\color{black}{$H_N = \text{span} 
\{e_0,e_1,e_2,\cdots,e_N\}$, $N\in \mathbb{N}^+$, where $\{e_i\}_{i=0}^N$ are  the $N+1$ first eigenvectors of $A$ in the space $H$.}}  
Then we propose the fully discrete method given by
\begin{align}\label{eq:full-discretization}
X_m^{M,N}-X^{M,N}_{m-1}
+
k A(A X_m^{M,N}
+
P_N F(X_m^{M,N}))
=
P_N \Delta W_m,\; X_0^{M,N}=P_NX_0.
\end{align}
Here $\Delta W_m := W(t_m)-W(t_{m-1})$ and $P_N$ is the orthogonal projector onto the space $H_N$.
As implied by Theorem \ref{them:error-estimates-full-stochastic}, the resulting spatio-temporal approximation error
is measured as follows
\begin{align}\label{eq:convergenc-rate}
\|X(t_m)-X_m^{M,N}\|_{L^p(\Omega;H)}=O(\lambda_N^{-\frac\gamma2}
+k^{\frac\gamma4}),\;
\gamma \in \big(\tfrac d 2, 4 \big],
\end{align}
{\color{black}{where $\lambda_N$ is the $N$-th eigenvalue of $A$ with the corresponding eigenvector $e_N$ .}}
This  reveals how the convergence rate depends on the regularity of the mild solution.  Comparing this with the spatial regularity
results \eqref{eq:sptial-regulaity}, one can clearly see that
 the order of convergence in space  coincides with the
spatial  regularity of the mild solution, for all $\gamma\in(\frac d2,4]$. However,  the convergence rate in time is more involved.
For $\gamma\in(\frac d2,2]$, the convergence rate in time agrees with the temporal H\"{o}lder regularity of the mild solution. For $\gamma\in[2,4]$,  the order of  the  convergence in time
 is higher than the H\"{o}lder regularity in time of the mild solution, this is due to the fact that that noise is additive noise. {\color{black}{This standard phenomenon is known for  SPDEs, see\cite{wang2016strong,kloeden2011exponential}.}}


Finally, we give some comments on a few closely relevant works about numerical methods  for the CHC equation.
  In \cite{kovacs2014erratum, larsson2011finite,furihata2018strong},
  strong convergence of mixed finite element methods for \eqref{eq:abstract-CHC-form} was proved, but with no rate obtained.
   The analysis in \cite{kovacs2014erratum, larsson2011finite,furihata2018strong} is based on proving a priori moment bounds with large exponents and
   in higher order norms using energy arguments and bootstrapping arguments followed by a pathwise Gronwall argument in the mild solution setting.
 Recently, the work \cite{ruisheng2018error}  recovers strong convergence rates of the mixed finite element methods (FEM)
 for the CHC equation with
spatial smooth noise, by using a priori moment bounds of the numerical approximations.
In \cite{cui2021strong}, the authors obtained convergence rates of  spectral Galerkin fully discrete schemes for
the CHC equation with space-time white noise.  {\color{black}{The paper \cite{cui2020absolute} presented 
the analysis of strong convergence rate of an implicit full discretization applied to stochastic Cahn-Hilliard equation  
with unbounded noise diffusiton in dimension one. The authors of \cite{feng2020fully} derived strong convergence rates of
a fully discrete mixed FEM for the stochastic Cahn--Hilliard equation with gradient-type multiplicative noise,
where the noise process is a real-valued Wiener process.}}

The outline of this paper is as follows. In the next two sections, some preliminaries are collected and well-posedness of the considered problem is elaborated. Section 4 is devoted to the uniform moment bounds of the fully discrete  approximation.
Based on the uniform moment bounds obtained in section 4, we derive the error estimates for the fully discrete problem in section 5.
%

\section{Settings}

In this section, we make some assumptions for the abstract equation \eqref{eq:abstract-CHC-form},
concerning the linear operator $A$,
the nonlinear term $F$, the noise process $W(t)$ and the initial data $X_0$.

\begin{assumption}\label{assum:linear-operator-A} (Linear operator $A$)
Let $D$ be a bounded convex domain in $\mathbb{R}^d$ for $d\in\{1,2,3\}$ with sufficiently smooth boundary and
let $H=L^2(D;\mathbb{R})$  be the real separable Hilbert space endowed
 with the usual inner product
 $\big<\cdot,\cdot\big>$
 and the associated norm
 $\|\cdot\|=\big<\cdot,\cdot\big>^{\frac12}$.
 Let
 $\dot{H}=\{v\in H:\int_D v\,\dd x=0\}$
 and
 $-A \colon dom(A)\subset \dot{H}\rightarrow \dot{H}$
 be  the Laplacian with homogeneous Neumann boundary conditions,
 defined by
$-Au = \Delta u$
 with
$
u\in dom(A) :=\{v\in H^2(D)\cap \dot{H}:\frac{\partial v}{\partial n}=0\}
$.
\end{assumption}

Assumption \ref{assum:linear-operator-A} guarantees that there exists a family of positive sequence
$\{\lambda_j\}_{j \in \mathbb{N}}$ and an orthonormal basis $\{e_j\}_{j \in\mathbb{N}}$ of $\dot{H}$
such that
$0<\lambda_1\leq\lambda_2\leq\cdots$,
$\lambda_j\sim j^{\frac2d} \rightarrow\infty $
and
$A e_j=\lambda_je_j$ for $j\in \mathbb{N}$. 
{\color{black}{Define $P:H\rightarrow \dot{H}$ as the orthogonal projector
such that
\begin{equation}
(I-P)v=|D|^{-1}\int_Dv\,\dd x.
\end{equation}
}}
When extended to $H$ as  $A v : = A P v$, for $v \in  H$,
the linear operator $A$ has an orthonormal basis $\{e_j\}_{j\in \mathbb{N}_0}$ of $H$  by taking $e_0=|D|^{-\frac12}$.
Throughout this paper,  we use $\N$ to denote the set of all positive integers and denote $\N_0 = \{ 0 \} \cup \N$.
Further, we introduce
\begin{align}
\dot{H}^\alpha=\big\{v\in \dot{H}: |v|_\alpha<\infty\big\},
\quad
H^\alpha=\big\{v\in H: \|v\|_\alpha<\infty\big\},
\end{align}
where the seminorm and norm are defined by
\begin{align}
|v|_\alpha
:=\big(
\sum_{j=1}^\infty \lambda_j^\alpha
|\big<
v,e_j
\big>|^2
\big)^{\frac12},
\;
\|v\|_\alpha
:=\big(
|v|_\alpha^2
+
|\big<v,e_0\big>|^2
\big)^{\frac12}.
\end{align}
By the spectral theory, we can also define the fractional powers
of $A$  on $\dot{H}$ in a simple way, e.g.,
$A^\alpha v=\sum_{j=1}^\infty\lambda_j^\alpha\left<v,e_j\right>e_j$, $\alpha\in \mathbb{R}$.
Note that
 $\dot{H}^\alpha=D(A^{\frac\alpha2})$
 and it is a real Hilbert space  with the inner
 product
 $\langle A^{\frac\alpha 2}\cdot, A^{\frac\alpha 2}\cdot \rangle$
 and the associated norm
 $
 \|\cdot\|_{\dot{H}^{\alpha}}
 :=
 |\cdot|_\alpha
 :=
 \|A^{\frac\alpha 2}\cdot\|
 $.
 It is well-known that for  integer $\alpha\geq 0$, the norm $\|\cdot\|_\alpha$ is  equivalent on $H^\alpha$ to the standard Sobolev norm $\|\cdot\|_{H^\alpha(D)}$, if $D$ is a regular domain in $\mathbb{R}^d$ with appropriately smooth boundary.

{\color{black}{In addition, Assumption \ref{assum:linear-operator-A}  ensures that  the operator $-A^2$ can
generate an analytic semigroup $E(t)=e^{-tA^2}$ on $H$, given by
\begin{align}
\begin{split}
E(t)v=e^{-tA^2}v
&=\sum_{j=0}^\infty e^{-t \lambda_j^2}\left<v,e_j\right>e_j
=
\sum_{j=1}^\infty e^{-t \lambda_j^2}\left<v,e_j\right>e_j
+
\left<v,e_0\right>e_0
\\
&=
{\color{black}{
e^{-tA^2}Pv}}
+
(I-P)v,
\quad
v \in H.
\end{split}
\end{align}}}
By expansion in terms of the eigenbasis of $A$ and using Parseval's identity,  one can easily obtain
\begin{align}
\| A^\mu E(t)\|_{\mathcal{L}(H)}
&
\leq
 Ct^{-\frac\mu2},\; t>0,\; \mu\geq 0,
 \label{I-spatio-temporal-S(t)}
 \\
\|A^{-\nu}(I-E(t))\|_{ \mathcal{L} (H)}
 &
 \leq
 Ct^{\frac\nu2},\quad t\geq0,\;\nu\in[0,2],
 \label{II-spatio-temporal-S(t)}
 \\
\int_{\tau_1}^{\tau_2} \|A^\varrho E(s) v\|^2\,\dd s
&
\leq
C|\tau_2-\tau_1|^{1-\varrho}\|v\|^2,\;\forall v\in H, \varrho\in[0,1],
\label{III-spatio-temporal-S(t)}
\\
 \Big\|A^{2\rho}\int_{\tau_1}^{\tau_2}
         E(\tau_2-\sigma)v\,\dd \sigma\Big\|
 &
 \leq
 C|\tau_2-\tau_1|^{1-\rho}\|v\|,\;\forall v\in H,\; \rho\in[0,1].
 \label{IV-spatio-temporal-S(t)}
 \end{align}
\begin{assumption}\label{assum:nonlinearity}(Nonlinearity) Let $F:L^6(D;\mathbb{R})\rightarrow H$ be a deterministic mapping given by
\begin{align}
F(v)(x)=f(v(x))=v^3(x)-v(x),
\quad
x\in D, v\in L^6(D;\mathbb{R}).
\end{align}
\end{assumption}
Here and below, $L^r(D;\mathbb{R}), r\geq 1$ ($L^r(D)$ or $L^r$ for short) is the Banach space
consisting of $r$-times integrable functions. Denote by $ V:=C(D;\mathbb{R}) $  the Banach space of continuous functions
with a usual norm. It is easy to check that, for any $v, \psi, \psi_1, \psi_2\in L^6(D)$,
\begin{align}\label{eq:definition-F-derivation}
\begin{split}
(F'(v)(\psi)\big)(x)&=f'(v(x))\psi(x)=\big(3v^2(x)-1\big)\psi(x),
\quad
x\in D,
\\
\big(F''(v)(\psi_1,\psi_2)\big)(x)&=f''(v(x))\psi_1(x)\psi_2(x)=6v(x)
\psi_1(x)\psi_2(x),
\quad
x\in D.
\end{split}
\end{align}
Moreover, one can check {\color{black}{the following coercivity, one-side Lipschitz  and polynomial growth conditions,}}
\begin{align}
\qquad\qquad
{\color{black}{-\big<F(u),u\big>}}
\leq
&
-\|u\|_{L^4}^4+\|u\|^2,\qquad u\in L^4,
\label{eq:cercivity-condition}
\\
-\big<F(u)-F(v),u-v\big>&\leq
\|u-v\|^2,\qquad\quad u, v\in L^6(D),
\label{eq:one-side-condition}
\\
\|F(u)-F(v)\|
\leq
\|u-v\| ( 1 + &\tfrac32\|u\|^2_V + \tfrac32\|v\|^2_V ),
\qquad
u, v\in V.
\label{eq:local-condition}
\end{align}

In order to define the $Q$-Wiener process, we introduce additional notations and spaces.
Let
$\mathcal{L}(H)$ be the Banach space of all bounded linear  operators from $H$ to $H$ endowed with the usual operator norm. Also, let $\mathcal{L}_2(H)$
be the Hilbert space consisting of  all Hilbert-Schmidt operators from
$H$ into $H$, equipped with the inner product and the norm,
\begin{align}
\big<\Gamma_1,\Gamma_2\big>_{\mathcal{L}_2(H)}
=
\sum_{j=1}^\infty \big<\Gamma_1 \phi_j,\Gamma_2\phi_j\big>,
\qquad
\|\Gamma\|^2_{\mathcal{L}_2(H)}
=
\sum_{j=1}^\infty \big\|\Gamma \phi_j\|^2,
\end{align}
where $\{\phi_j\}_{j\in\N}$ is an arbitrary orthonormal basis  of $H$.
If $\Gamma\in\mathcal{ L}_2(H)$ and $L\in \mathcal{L}(H)$, then $\Gamma L, L\Gamma\in \mathcal{L}_2(H)$ and
\begin{align}
\|\Gamma L\|_{\mathcal{L}_2(H)}
\leq
\|L\|_{\mathcal{L}(H)}
\|\Gamma\|_{\mathcal{L}_2(H)}
,
\|L\Gamma \|_{\mathcal{L}_2(H)}
\leq
\|L\|_{\mathcal{L}(H)}
\|\Gamma\|_{\mathcal{L}_2(H)}.
\end{align}
\begin{assumption}\label{assum:noise-term}(Noise process)
Let $\{W(t)\}_{t\in[0,T]}$ be a standard $H$-valued $Q$-Wiener process with the covariance operator
 $ \mathcal{L}(H)  \ni Q:H \rightarrow H $ being a  symmetric nonnegative operator satisfying
{\color{black}{\begin{align}\label{eq:ass-AQ-condition}
\|A^{\frac{\gamma - 2}2}PQ^{\frac12}\|_{\mathcal{L}_2(H)}
<\infty,
\quad
\text{ for some } \quad \gamma\in \big(\tfrac d2,4\big].
\end{align}
}}
\end{assumption}
{\color{black}{
We mention that the assumption of  $\gamma > \tfrac d2$ is necessary in the following analysis. 
For the space-time white noise $(Q = I)$, the condition  \eqref{eq:ass-AQ-condition} is satisfied
for $\gamma < 2 - \tfrac{d}{2} $,  see \cite{larsson2011finite}.
Since we also demand $\gamma > \tfrac d2$, the study of the space-time white noise case is thus
limited to dimension $1$. For the trace-class noise ($\text{Tr} (Q) = \| Q^{\frac12}\|_{\mathcal{L}_2(H)}< \infty$), 
the assumption \eqref{eq:ass-AQ-condition} is fulfilled with $\gamma = 2$ and multi-dimensions $d \geq 2$ are all allowed.
}}
\begin{assumption}\label{assum:intial-value-data}(Initial data)
Let $X_0:\Omega \rightarrow H$ be $\mathcal{F}_0/\mathcal{B}(H)$-measurable and satisfy, for a sufficiently large number $p_0\in \mathbb{N}$,
\begin{align}
{\color{black}{\mathbb{E}[\|X_0\|_\gamma^{p_0}]}}<\infty,
\end{align}
where $\gamma\in\big( \tfrac d2,4\big] $ is the parameter from \eqref{eq:ass-AQ-condition}.
\end{assumption}

\section{Well-posedness and regularity of the CHC equation}
%
The aim of this section is to elaborate the well-posedness and
spatio-temporal regularity of the mild solution to the considered problem.
%
To this end,
we split the mild solution as $X(t)=Y(t)+Z(t)$,
 where
 $Z(t)$
 satisfies
 {\color{black}{\begin{align}\label{eq:stochastic-convulution-problem}
 \,\dd Z(t)+A^2Z(t)\,\dd t=\dd W(t),&\; Z(0)=(I-P)X_0,
 \end{align}}}
and $Y(t)$ satisfies
  \begin{align}\label{eq:abstract-CHC-form-I}
 \dot{Y}(t)
+
A^2Y(t)+APF(Y(t)+Z(t))
=
0,\; Y(0)={\color{black}{PX_0}}.
\end{align}
 It is well-known that the problem \eqref{eq:stochastic-convulution-problem} admits a unique mild solution, given by
 \begin{align}\label{mild-solution-SCP}
{\color{black}{Z(t)=E(t)(I-P)X_0+ \mathcal{O}(t)}}
,
 \end{align}
{\color{black}{where $\mathcal{O}(t):=\int_0^t  E(t-s)  \mathrm{d} W(s)$ is called the stochastic convolution.}}
 In the next lemma, we consider the regularity of $Z(t)$. For the proof of spatial regularity result \eqref{lemma:p-th-moment-bound-mild-solution}, we refer to \cite[Theorem 3.1]{larsson2011finite}.
 By the similar arguments of the proof of  \cite[Theorem 2]{chai2018conforming}, we also obtain the following temporal regularity result.
 \begin{lemma}\label{lemma:reguralirty-mild-convolution}
If Assumptions \ref{assum:linear-operator-A},\ref{assum:noise-term} are valid,  then {\color{black}{the mild solution of the problem \eqref{eq:stochastic-convulution-problem}, given by \eqref{mild-solution-SCP},}}  satisfies the following spatio-temporal regularity results as follows, for $p\geq 1$,
\begin{align}\label{lemma:p-th-moment-bound-mild-solution}
\mathbb{E}
{\color{black}{\Big[
\sup_{s\in[0,T]}\|Z(s)\|_\gamma^p
\Big]
<
\infty,}}
\end{align}
and for $\beta\in[0,\gamma]$,
\begin{align}\label{lemma:temporal-stochastic-conlution}
{\color{black}{\|Z(t)-Z(s)\|_{L^p(\Omega;H^\beta)}
\leq
C|t-s|^{ \min \{\frac12,\frac{\gamma-\beta}4\} }.}}
\end{align}
\end{lemma}
 Therefore, it suffices to treat  \eqref{eq:abstract-CHC-form-I}. In the next two subsections, we will use Galerkin's method and energy arguments to address this issue.

 \subsection{Useful inequalities}
In this part, we collect  some useful inequalities, which play an important role in our error analysis below. We first introduce the following  embedding inequalities,
\begin{align}\label{eq:embedding-equatlity-I}
{\color{black}{H^1\subset H^{\frac d3}}}\subset L^6(D)\quad \text{and} {\color{black}{\quad H^\delta}} \subset C(D;\mathbb{R}),
\quad
\text{ for } \delta>\tfrac d2, \;
d\in\{1,2,3\}.
\end{align}
{\color{black}{In light of the definition of $P$}}, it is not difficult to check
\begin{align}
\|Pv\|_{L^q}
\leq 
{\color{black} (1+|D|^{-\frac1q}) }
\|v\|_{L^q},\;q\geq 2.
\end{align}
This together with \eqref{eq:embedding-equatlity-I} implies  (see \cite{ruisheng2018error}), for any $\delta>\frac d2$, 
\begin{align}\label{eq:embedding-equatlity-III}
\begin{split}
\|A^{-\frac\delta2}P x\|
\leq
C\|x\|_{L^1},
\;
\|A^{-\frac12} Px\|
\leq
C\|x\|_{L^{\frac65}},\;\forall x\in L^2(D).
\end{split}
\end{align}
Additionally, we also use the fact that,
 the norm $\|\cdot\|_2$ is  equivalent on $H^2$ to
the standard Sobolev norm $\|\cdot\|_{H^2(D)}$ {\color{black}{(see \cite{ee2006galerkin, yagi2009abstract})}} to derive, for any $f,g\in H^2$,
\begin{align}\label{eq:algebra-properties-Hs}
\|fg\|_{H^2(D)}
\leq
C\|f\|_{H^2(D)}\|g\|_{H^2(D)}
\leq
C{\color{black}{\|f\|_2\,\|g\|_2}}.
\end{align}

Based on the above preparations, we derive the following inequality.
\begin{lemma}\label{lem:Fx1-Fx2}
Let $g:L^4\rightarrow H$ be the Nemytskij operator of a polynomial of second degree. The following estimates hold
 for any {\color{black}{$\iota\in (\frac12,1)$}} and $d=1$,
\begin{align}\label{eq:FX1-FX2-d=1}
\|g(u)v\|_{{\color{black}{\iota}}}
\leq
C\|v\|_1(1+\|u\|_{\iota}^2), \text{ for } \;u\in {\color{black}{H^\iota}}, v\in H^1,
\end{align}
and
for any {\color{black}{$\iota\in(\frac d2,2)$}},  $d=2,3$,
\begin{align}\label{eq:FX1-FX2-d=23}
\|g(u)v\|_{1}
\leq
C\|v\|_{1}(1+\|u\|_{\iota}^2),
for\;u\in H^\iota, v\in H^1.
\end{align}
\end{lemma}
{\it Proof of Lemma \ref{lem:Fx1-Fx2}.}
First, we consider the case $d=1$. The H\"{o}lder inequality,  \eqref{eq:embedding-equatlity-I} and {\color{black}{the equivalence of norms in $H^\iota$ and $H^\iota(D)$}} yield that for  $u, v\in H^\iota$,
\begin{align}\label{eq:lem:Fx1-Fx2-one-dimen}
\begin{split}
\|g(u)v\|^2_{\iota}
&\leq
\|g(u)v\|^2
+
C\int_D\int_D\frac{|g(u(\xi_1))v(\xi_1)-g(u(\xi_2))v(\xi_2)|^2}
{|\xi_1-\xi_2|^{2\iota+1}}\,\dd \xi_1\,\dd \xi_2
\\
&
\leq
C\|v\|^2(1+\|u\|_V^4)
+
C(1+\|u\|_V^4)\|v\|^2_{\iota}
+
C\|u\|^2_{\iota}(1+\|u\|_V^2)\|v\|^2_V
\\
&
\leq
C(1+\|u\|_V^4+\|u\|^4_{\iota})(\|v\|^2_V+\|v\|_{\iota}^2)
\\
&
\leq
C(1+\|u\|^4_{\iota})\|v\|_{\iota}^2.
\end{split}
\end{align}
 The above estimate assures \eqref{eq:FX1-FX2-d=1}.
For the case $d=2,3$, we first obtain
\begin{align}\label{eq:lem:Fx1-Fx2-two-three-dim}
\begin{split}
\|g(u)v\|_{1}
&\leq
C(\|g(u)v\|
+
\|v \nabla g(u) \|
+
\|g(u) \nabla v\|)
\\
&\leq
C(1+\|u\|^2_V)\|v\|
+
\| v \nabla g(u)\|
+
C\|v\|_1(1+\|u\|_V^2)
.
\end{split}
\end{align}
Further,
owing to \eqref{eq:embedding-equatlity-I} and the fact that the Sobolev embedding theorem tells us  ${\color{black}{H^s}}\subset L^p$, for $0<s<1$ and $p=\frac{2d}{d-2s}$,
there exists a constant $\delta$ given by
$
{\color{black}{\delta=\frac{2\iota-2}{2-\iota}}}$
for $d=2$ and $\delta=1$ for $d=3$ such that, for $\iota>\frac d 2, d=2,3$
{\color{black}{
\begin{align}\label{eq:inequality-v-u}
\|v\nabla u\|
\leq
C\|v\|_{L^{\frac {2(2+\delta)} \delta}} \|\nabla u\|_{L^{2+\delta}}
\leq
C\|v\|_1 \| u\|_{H^\iota}
\end{align}}}
Therefore,
\begin{align}\label{eq:bound-F2-p1}
\begin{split}
\|v\nabla g(u) \|
&\leq
C\|v\nabla u \|(1+\|u\|_V )
\leq
C\|v\|_1 (1+\|u\|_{\iota}^2 ).
\end{split}
\end{align}
This together with \eqref{eq:embedding-equatlity-I} and \eqref{eq:lem:Fx1-Fx2-two-three-dim} can show \eqref{eq:FX1-FX2-d=23}. Hence this ends the proof.
$\square$

\subsection{The semidiscrete Galerkin approximation}
In this part, we use the Galerkin method to approximate the problem \eqref{eq:abstract-CHC-form-I}. For $n\in \mathbb{N}$,  we define a finite dimensional subspace of ${\color{black}{H}}$ by
\begin{align}
H_n=
{\color{black}{span\{e_0,e_1,e_2\cdots,e_n\}}}
\end{align}
and the projection operator $P_n: {\color{black}{H^\alpha}}\rightarrow H_n$ by $P_n\xi={\color{black}{\sum_{i=0}^{n}\big<\xi,e_i\big>e_i}}$, for $\forall \xi\in H^\alpha$, $\alpha\geq -2$. Here $\{e_i\}_{i=0}^{n}$ are  the {\color{black}{$n+1$}} first eigenvectors of the dominant linear operator $A$ in the space $H$. It is not difficult to find that the operators $A$ and $P_n$ {\color{black}{commute, but $P_nP\neq P$}}. Moreover,
\begin{align}\label{eq:property-Pm}
\|(I-P_n)\varphi\|
\leq
\lambda_{n+1}^{-\frac\alpha2}|\varphi|_\alpha, \;\forall \varphi\in H^\alpha, \;\alpha>0.
\end{align}
Then
the Galerkin approximation of \eqref{eq:abstract-CHC-form-I} is  given by
\begin{align}\label{eq:gerlrkin-approximation-CHC-i}
\begin{split}
\left\{
\begin{array}{ll}
\frac{\dd v^{n}}{\dd t}+A^2v^{n} +P_nAP F(v^{n}+{\color{black}{Z^{n}}})=0,&\; \text{ in }\; D\times (0,T],
\\
v^{n}(0)={\color{black}{P_nPX_0}},&\; \text{ in } \; D,
\end{array}
\right.
\end{split}
\end{align}
where ${\color{black}{Z^{n}(t)=P_nZ(t)}}$.
It is clear that the problem \eqref{eq:gerlrkin-approximation-CHC-i} admits a unique solution
given by
\begin{align}\label{eq:mild-solution-gerlrkin-approximation-CHC-i}
v^{n}(t)
=
E(t)P_nPX_0
-
\int_0^tE(t-s)P_nAPF(v^{n}(s)+{\color{black}{Z^n(s)}})\,\dd s,
\quad
t \in [0, T].
\end{align}
Here $Z^{n}(t)$ enjoys the similar  spatio-temporal regularity as  $Z(t)$:
\begin{align}\label{lemma:p-th-moment-bound-mild-galerkin-solution}
\sup_{n\in \mathbb{N}}
\mathbb{E}
\Big[
\sup_{s\in[0,T]}
\big\|
{\color{black}{Z^{n}(t)}}
\big\|_\gamma^p
\Big]
<
\infty,
\end{align}
and for $\beta\in[0,\gamma]$,
\begin{align}\label{lemma:temporal-mild-galerkin-solution}
\sup_{n\in \mathbb{N}}
\|{\color{black}{Z^{n}(t)-Z^{n}(s)}}\|_{L^p(\Omega;H^\beta)}
\leq
C|t-s|^{\min\{\frac12,\frac{\gamma-\beta}4\}}.
\end{align}
  {\color{black}{Setting}} $X^{n}(t):=v^{n}(t)+ {\color{black}{Z^{n}}}(t)$, we will show that the sequence {\color{black}{$\{X^{n}(\cdot)\}$}} is Cauchy in $C_W([0,T];H)$, where $C_W([0,T];H)$ is a Banach space consisting of all continuous mappings $G \colon [0,T]\rightarrow L^2(\Omega,\mathcal{F}, \mathbb{P};H)$, endowed with the norm
 \begin{align}
 \| G \|_{C_W([0,T];H)}=
 \Big(
 \sup_{s\in[0,T]}\mathbb{E}\big[ \| G ( s ) \|^2\big]
 \Big)^{\frac12}.
  \end{align}
{\color{black} 
To prove that the limit of $\{X^{n}(t)\}$ in $C_W([0,T];H)$ is the required mild solution of \eqref{eq:abstract-CHC-form},
let us start with {\color{black}{a priori estimate}} of $X^n(t)$.
}
 \begin{lemma}\label{eq:approximate-solution-H}
Let $v^n(t)$ be a solution of \eqref{eq:gerlrkin-approximation-CHC-i} and let $X^n(t) := v^n(t)+ {\color{black}{Z^{n}(t)}}$.
If Assumptions \ref{assum:linear-operator-A}-\ref{assum:intial-value-data} are valid, then  there exists a constant $C=C(T,\gamma,p)$ such that
 \begin{align}\label{eq:projection-property-projection-galkerin}
 \sup_{n\in \mathbb{N}}
 \mathbb{E}
 \Big[
 \sup_{s\in[0,T]}\|X^{n}(s)\|^p
 \Big]
 \leq C
 <\infty.
 \end{align}
 \end{lemma}
{\color{black}Its proof is not new and has been already used in the literature. We keep it here for completeness.}

 {\it Proof of Lemma \ref{eq:approximate-solution-H}.}
Multiplying both sides of \eqref{eq:gerlrkin-approximation-CHC-i} by $A^{-1}v^n$ and {\color{black}{using \eqref{eq:cercivity-condition}}} and the chain rule
and the fact 
 yield
\begin{align}\label{eq:vm-H-1-norm}
\begin{split}
|v^{n}(t)|_{-1}^2
\leq
&
|{\color{black}{P_nPX_0}}|_{-1}^2
-
2\int_0^t
|v^{n}(s)|_1^2\,\dd s
-2
\int_0^t\big<   F(v^{n}(s)+{\color{black}{Z^{n}(s)}},v^{n}(s)\big>\,\dd s
\\
\leq
&
|X_0|_{-1}^2
-
2\int_0^t
|v^{n}(s)|_1^2\,\dd s
-
2\int_0^t\|v^{n}(s)\|^4_{L^4}\,\dd s
+
2\int_0^t\|v^{n}(s)\|^2_{L^2}\,\dd s
\\
&-
2\int_0^t{\color{black}{\big<3(v^{n}(s))^2Z^{n}(s)+3(Z^{n}(s))^2v_n(s)
+(Z^{n}(s))^3-Z^{n}(s),v^{n}(s)\big>}}\,\dd s
\\
\leq
&
|X_0|_{-1}^2
-
2\int_0^t
|v^{n}(s)|_1^2\,\dd s
-
2\int_0^t\|v^{n}(s)\|^4_{L^4}\,\dd s
+
2\int_0^t\|v^{n}(s)\|^2_{L^2}\,\dd s
\\
&
{\color{black}{+
6\int_0^t\|v^n(s)\|_{L^4}^3\|Z^n(s)\|_{L^4}\,\dd s
+
6\int_0^t\|v^n(s)\|_{L^4}\|Z^n(s)\|_{L^4}^2\,\dd s}}
\\
&
{\color{black}{+
\int_0^t\|v^n(s)\|_{L^4}^2\|Z^n(s)\|_{L^4}^2\,\dd s
+
\int_0^t\|v^n(s)\|\|Z^n(s)\|\,\dd s}}
\\
\leq
&
|X_0|_{-1}^2
-
2\int_0^t
|v^{n}(s)|_1^2\,\dd s
-
\int_0^t\|v^{n}(s)\|^4_{L^4}\,\dd s
+
C\int_0^t(1+\|{\color{black}{Z^{n}(s)}}\|_{L^4}^4)\,\dd s,
\end{split}
\end{align}
{\color{black}{where Young's inequality was also used  in the last inequality.}}
Therefore,
\begin{align}
2\int_0^t|v^{n}(s)|_1^2\,\dd s
+
\int_0^t\|v^{n}(s)\|_{L^4}^4\,\dd s
\leq
|X_0|^2_{-1}
+
C\int_0^t(1+{\color{black}{\|Z^{n}(s)\|_{L^4}^4}})\,\dd s,
\end{align}
which in combination with  \eqref{eq:embedding-equatlity-I} and  \eqref{lemma:p-th-moment-bound-mild-galerkin-solution} leads to
\begin{align}\label{eq:L4-norm-integrand}
\begin{split}
\mathbb{E}
\Big[
\int_0^t
|v^{n}(s)|_1^2
\,\dd s
\Big]^p
+
\mathbb{E}
\Big[
\int_0^t
\|v^{n}(s)\|_{L^4}^4
\,\dd s
\Big]^p
&\leq
C
\Big(
\mathbb{E}
\big[
|X_0|_{-1}^{2p}
\big]
+
\int_0^t
\big(1+\mathbb{E}
\big[
\|{\color{black}{Z^{n}(s)}}\|_V^{4p}
\big]
\big)\;\dd s
\Big)
\\
&\leq
C
\Big(
\mathbb{E}
\big[
|X_0|_{-1}^{2p}
\big]
+
\int_0^t
\big(
1+\sup_{n\in \mathbb{N}}
\mathbb{E}
\big[
{\color{black}{\|Z^{n}(s)\|_\gamma^{4p}}}
\big]
\big)\;\dd s
\Big)
\\
&
<\infty.
\end{split}
\end{align}
Now we are in the position to give the a priori estimate of $\|v^n\|$.
Taking inner product of \eqref{eq:gerlrkin-approximation-CHC-i} by $v^n(t)$ and using the chain rule and {\color{black}{the self-adjointness of $A$}}, one can find that
\begin{align}
\begin{split}
\| v^{n}(t)\|^2
=
&
\|P_nP X_0\|^2
-
2 \int_0^t  \|A v^{n}(s)\|^2  \,\dd s
+
2 \int_0^t
{\color{black}{
\big<
F
\big(
v^{n}(s)
+
Z^n(s)
\big)
,
-A v^n(s)
\big>
\,\dd s
}}
\\
=
&
{\color{black}{
\|P_nP X_0\|^2
-
2 \int_0^t  \|A v^{n}(s)\|^2  \,\dd s
+
2 \int_0^t
\big<
\big(
v^{n}(s))^3-v^{n}(s)
\big)
,
-A v^n(s)
\big>
\,\dd s
}}
\\
&
+
2
\int_0^t
\big<
{\color{black}{3(v^{n}(s))^2Z^{n}(s)
+
3v^{n}(s)(Z^{n}(s))^2
+
(Z^{n}(s))^3
-
Z^{n}(s)}}
,
{\color{black}{-Av^{n}(s)}}
\big>
\,\dd s
\\
\leq
&
\|X_0\|^2
-
2 \int_0^t
\|Av^{n}(s)\|^2
\,\dd s
-
\tfrac32
\int_0^t
|(v^{n}(s))^2|_1^2
\,\dd s
+
2
\int_0^t
| v^{n}(s)|_1^2
\,\dd s
\\
&+
{\color{black}{
2
\int_0^t
\big(
\|{\color{black}{Z^{n}(s)}}\|_V^2
\|v^{n}(s)\|^2_{L^4}
+
\| v^{n}(s)\|
\|{\color{black}{Z^{n}(s)}}\|_V^2
+
\|{\color{black}{Z^{n}(s)}}\|_{L^6}^3
+
\|{\color{black}{Z^{n}(s)}}\|)\|Av\|
\,\dd s
}}
\\
\leq
&
\|X_0\|^2
-
\int_0^t
\|Av^{n}(s)\|^2
\,\dd s
-
\tfrac32
\int_0^t
|(v^{n}(s))^2|_1^2
\,\dd s
+
2
\int_0^t
| v^{n}(s)|_1^2
\,\dd s
\\
&+
C
\int_0^t
\big(
\|{\color{black}{Z^{n}(s)}}\|_V^4
\|v^{n}(s)\|^4_{L^4}
+
\| v^{n}(s)\|
\|{\color{black}{Z^{n}(s)}}\|_V^4
+
\|{\color{black}{Z^{n}(s)}}\|_{L^6}^6
+
\|{\color{black}{Z^{n}(s)}}\|^2)
\,\dd s.
\end{split}
\end{align}
Hence, we have
\begin{align}
\begin{split}
&
\|v^{n}(t)\|^2
+
\int_0^t \!\!
\|Av^{n}(s)\|^2
\,\dd s
+
\tfrac32
\int_0^t \!\!
|(v^{n}(s))^2|_1^2
\,\dd s
\\
&
\leq
C
\|X_0\|^2
+
C
\int_0^t\!\!
\big(
1+\|{\color{black}{Z^{n}(s)}}\|_V^6
\big)
\big(
1+\|v^{n}(s)\|^4_{L^4}
\big)
\,\dd s.
\end{split}
\end{align}
Further, from  \eqref{eq:L4-norm-integrand}, \eqref{lemma:p-th-moment-bound-mild-galerkin-solution} and  the H\"{o}lder inequality,
it follows that for any $n\in \mathbb{N}$,
\begin{align}\label{eq:L2-Y-and-H2-Y-integrand}
\begin{split}
\mathbb{E}
\Big[
\sup_{s\in[0,T]}
&
\|v^{n}(s)\|^{2p}
\Big]
+
\mathbb{E}
\Big[
\int_0^T
\|Av^{n}(s)\|^2
\,\dd s
\Big]^p
+
\mathbb{E}
\Big[
\int_0^T
|(v^{n}(s))^2|_1^2
\,\dd s
\Big]^p
\\
\leq
&
C(p)~
\mathbb{E}
[\|X_0\|^{2p}]
+
C(p,T)
\Big(
1
+
\mathbb{E}
\Big[
\sup_{s\in[0,T]}
\|{\color{black}{Z^{n}(s)}}
\|_\gamma^{12p}
\Big]
\Big)^{\frac12}
\Big (
1+
\mathbb{E}
\Big[
\int_0^T
\|v^{n}(s)\|^4_{L^4}
\,\dd s
\Big]^{2p}
\Big)^{\frac12}
\\
\leq
&C(p,T,X_0).
\end{split}
\end{align}
This together with \eqref{lemma:p-th-moment-bound-mild-galerkin-solution} finishes the proof of this lemma. $\square$
%

{\color{black}{Moreover, we  provide a  stronger moment bound and a   temporal regularity result of $X^n(t)$, 
which will play an important role in the proof of the existence and   uniqueness  of the mild solution $X(t)$.}}
\begin{lemma}\label{eq:approximate-solution-higher-norm}
Let $v^n(t)$ be the solution of \eqref{eq:gerlrkin-approximation-CHC-i}
and
$X^n(t)=v^n(t)+ {\color{black}{Z^{n}(t)}}$.
If Assumptions \ref{assum:linear-operator-A}-\ref{assum:intial-value-data} are valid, then we obtain
 \begin{align}\label{lemma:spatial-galerkin-solution}
\sup_{n\in \mathbb{N}}
\sup_{s\in[0,T]}
\|X^{n}(s)\|_{L^p(\Omega;H^\gamma)}
 <\infty,
 \end{align}
 and for any $\beta\in[0,\gamma]$,
  \begin{align}\label{lem:temporal-galerkin-solution}
 \sup_{n\in \mathbb{N}}
 \|X^{n}(t)-X^{n}(s)\|_{L^p(\Omega;H^\beta)}
 \leq
 C|t-s|^{\min\{\frac12,\frac{\gamma-\beta}4\}}.
 \end{align}
 \end{lemma}
 {\it Proof of Lemma \ref{eq:approximate-solution-higher-norm}.}
{\color{black}{The proof of this lemma in dimension $1$ can be given by 
following  the steps of the proof of \cite[Proposition 3.1]{cui2021strong}. 
In what follows we focus on the proof of \eqref{lemma:spatial-galerkin-solution}-\eqref{lem:temporal-galerkin-solution} 
in dimensions two and three.}}
For the  case
$\gamma\in [3,4]$,
please refer to \cite{ruisheng2018error}.
Therefore, it suffices to show the case
$\gamma\in (\frac d2,3]$.
By taking any fixed number
$\delta_0\in(\frac32,2)$,
we consider two possibilities:
either
$\gamma\in(\frac d2,\delta_0]$
or
$\gamma\in(\delta_0,3]$.
, we introduce a Lyapunov functional $J(u)$, defined by
\begin{align}\label{eq:definition-lyapunov-function}
J(u)=\frac12\|\nabla u\|^2
+
\int_D\Phi(u)\,\dd x,
\end{align}
where $\Phi$ is the primitive of $F$ vanishing at zeros, i.e., {\color{black}{$\Phi(s)=\int_0^sF(t) \,\dd t, s \in \R$}}.
We multiply \eqref{eq:gerlrkin-approximation-CHC-i} by $A^{-1}\dot{v}^n$ to arrive at
\begin{align}\label{eq:galerkin-modified-equation}
\begin{split}
\big|\dot{v}^{n}(t)\big|_{-1}^2
+
\tfrac12 \tfrac{\dd |v^{n}(t)|_1^2}{\dd t}
+
\big<F(v^{n}(t)+{\color{black}{Z^{n}(t)}}),\;\dot{v}^{n}(t)\big>
=
0,
\end{split}
\end{align}
{\color{black}
where $v^n( \cdot )$ is a solution of \eqref{eq:gerlrkin-approximation-CHC-i} and $\dot{v}^{n}( \cdot )$
stands for the time derivative of $v^n( \cdot )$.
}
Using the fact $\Phi'(u)=F(u)$ {\color{black}{and $P\dot{v}_n=\dot{v}_n$}} yields
\begin{align}
\begin{split}
&
-
\big<
F(v^{n}(t)+{\color{black}{Z^{n}(t))}}
,
\;
\dot{v}^n(t)
\big>
\\
=
&
-
\big<
F(v^{n}(t))
,\;
\dot{v}^n(t)
\big>
-
\big<
{\color{black}{3(v^{n}(t))^2Z^{n}(t)
+
3v^{n}(t)(Z^{n}(t))^2
+
(Z^{n}(t))^3-Z^{n}(t)}}
,\;
\dot{v}^{n}(t)
\big>
\\
\leq
&
-
\tfrac{\dd }{\dd t}\!
\int_D \!\!\!
\Phi(v^{n}(t))
\,\dd x
+
\big|
P
\big(
{\color{black}{3(v^{n}(t))^2Z^{n}(t)
+
3v^{n}(t)(Z^{n}(t))^2
+
(Z^{n}(t))^3
-
Z^{n}(t)}}
\big)
\big|_1
|\dot{v}^{n}(t)|_{-1}.
\end{split}
\end{align}
To treat the above {\color{black}{expression}}, we 
use \eqref{eq:embedding-equatlity-I} and \eqref{eq:inequality-v-u} with $\kappa=\gamma$ to obtain, for $\gamma>\frac d 2$, $d=2,3$,
\begin{align}\label{eq:v-convulution-in-h1-norm}
\begin{split}
\big|
P
\big(
3( v^{n})^2 {\color{black}{Z^{n}}}
&+
3 v^{n}{\color{black}{(Z^{n})^2}}
+
{\color{black}{( Z^{n})^3
-
Z^{n}}}
\big)
\big |_1
\\
&
{\color{black}{
\leq
\big\|
\nabla\big(
3( v^{n})^2 {\color{black}{Z^{n}}}
+
3 v^{n}{\color{black}{(Z^{n})^2}}
+
{\color{black}{( Z^{n})^3
-
Z^{n}}}
\big)
\big \|}}
\\&
\leq
C
\big(
\|(v^{n})^2 \|_1
\;
\|Z^{n}\|_V
+
\|(v^{n})^2\nabla Z^{n}\|
+
\|v^n\|_1
\;
\|Z^{n}\|_V^2
\\
&
\quad
+
\|v^{n}\|_V
\| Z^{n}\|_1
\;
\|Z^{n}\|_V
+
\| Z^{n}\|
\;
\|Z^{n}\|_V^2
+
\| Z^{n}\|_1
\big)
\\
&
{\color{black}{\leq
C
\big(
\|(v^{n})^2 \|_1
\;
\|Z^{n}\|_\gamma
+
\|(v^{n})^2\|_1\| Z^{n}\|_\gamma
+
\|v^n\|_1
\;
\|Z^{n}\|_\gamma^2}}
\\
&
{\color{black}{\quad
+
\|Av^{n}\|
\;
\|Z^{n}\|^2_\gamma
+
\| Z^{n}\|_\gamma^3
+
\| Z^{n}\|_\gamma
\big)}}
\\&
\leq
C(\gamma)
\big(
1+
\|(v^{n})^2 \|_1
+
\|A v^{n}\|
\big)
\big(
1+\|Z^{n}\|_\gamma^3
\big)
\end{split}
\end{align}
Therefore,
\begin{align}
\begin{split}
-\big<F(v^{n}(t)+Z^{n}),\;\dot{v}^{n}(t)\big>
\leq
&
-\frac{\dd }{\dd t}\int_D\Phi(v^{n}(t))\,\dd x
+
\frac12\|\dot{v}^{n}(t)\|_{-1}^2
\\
&
+
C(\gamma)\big(1+\|(v^n(t))^2\|_1+\|A v^{n}(t)\|\big)^2\big(1+\|Z^{n}(t)\|_\gamma^3\big)^2.
\end{split}
\end{align}
The above estimate together with \eqref{eq:galerkin-modified-equation} and \eqref{eq:definition-lyapunov-function} yields
\begin{align}
\begin{split}
J(v^{n}(t))
\leq
J({\color{black}{P_nPX_0}})
+
C(\gamma,T)
\Big(
1
+
\int_0^t \!
\big(
\|(v^{n}(s))^2 \|_1^2
+
\|A v^{n}(s)\|^2
\big)
\,\dd s
\Big)
\Big(
1+\sup_{s\in[0,T]}
\big\|
{\color{black}{Z^{n}(s)}}
\big\|_{\gamma}^{6}
\Big).
\end{split}
\end{align}
Then applying
 \eqref{eq:L2-Y-and-H2-Y-integrand}, \eqref{eq:L4-norm-integrand}, \eqref{eq:definition-lyapunov-function}
 and
 \eqref{lemma:p-th-moment-bound-mild-galerkin-solution}
 infers that
\begin{align}
\sup_{n\in \mathbb{N}}
\sup_{s\in[0,T]}
\|v^{n}(s)\|_{L^{2p}(\Omega;H^1)}
<
\infty,
\end{align}
which together with  \eqref{lemma:p-th-moment-bound-mild-solution} arrives at, for $d=2,3$,
\begin{align}\label{eq:bound-F-H}
\sup_{n\in \mathbb{N}}
\sup_{s\in[0,T]}
\|F(X^{n}(s))\|_{L^p(\Omega;H)}
<
\infty.
\end{align}
Hence,  when \eqref{eq:ass-AQ-condition} is fulfilled with $\gamma\in(\frac d2,\delta_0]$, we utilize
\eqref{I-spatio-temporal-S(t)} with $\mu=0$ and $1+\frac{\gamma}2$ to show
\begin{align}\label{eq:sptial-bound-vm-gamma-low-condition}
\begin{split}
\|v^{n}(t)\|_{L^p(\Omega;H^\gamma)}
&
\leq
\|E(t)P_nPX_0\|_{L^p(\Omega;H^\gamma)}
+
\Big\|
\int_0^t
E(t-s) P_n A P F(X^{n}(s))
\,\dd s
\Big\|_{L^p(\Omega;H^\gamma)}
\\
&
\leq
\|X_0\|_{L^p(\Omega;H^\gamma)}
+
C
\int_0^t
(t-s)^{-\frac{2+\gamma}4}
\,\dd s
\sup_{n\in \mathbb{N}}
\sup_{s\in[0,T]}
\|F(X^{n}(s))\|_{L^p(\Omega;H)}
\\
&
\leq
C(T,\gamma,X_0,p).
\end{split}
\end{align}
For the temporal regularity of  $v^{n}$, we apply \eqref{I-spatio-temporal-S(t)}, \eqref{II-spatio-temporal-S(t)} and
\eqref{eq:bound-F-H} to obtain, for $\beta\in [0,\gamma]$ and $t>s$,
\begin{align}\label{eq:temporal-bound-vm-gamma-low-condition}
\begin{split}
\|v^{n}(t)-v^{n}(s)\|_{L^p(\Omega;H^\beta)}
\leq
&
\|(E(t-s)-I)v^{n}(s)\|_{L^p(\Omega;H^\beta)}
+
\Big\|
\int_s^t
E(t-r) P_n A P F(X^{n}(r))
\,\dd r
\Big\|_{L^p(\Omega;H^\beta)}
\\
\leq
&
C(t-s)^{\frac{\gamma-\beta}4} \|v^{n}(s)\|_{L^p(\Omega;H^\gamma)}
+
\int_s^t
(t-r)^{-\frac{2+\beta}4}
\|
F(X^{n}(r))
\|_{L^p(\Omega;H)}
\,\dd r
\\
\leq
&
C(t-s)^{\frac{\gamma-\beta}4}
\sup_{n\in \mathbb{N}}
\sup_{s\in[0,T]}
\Big(\|v^{n}(s)\|_{L^p(\Omega;H^\gamma)}
+
\|F(X^{n}(s))\|_{L^p(\Omega;H)}
\Big)
\\
\leq
&C(t-s)^{\frac{\gamma-\beta}4}.
\end{split}
\end{align}
Combining the above two estimates with \eqref{lemma:p-th-moment-bound-mild-galerkin-solution}
and \eqref{lemma:temporal-mild-galerkin-solution} shows the case $\gamma\in (\frac d2,\delta_0]$.

Next, let us look at the other case $\gamma\in (\delta_0,3]$.
As already verified in the former case, one can see
$\sup_{n\in \mathbb{N}}
\sup_{s\in[0,T]}
\|X^n(s)\|_{L^p(\Omega;H^{\delta_0})}<\infty$
in this case.
Therefore, by utilizing  \eqref{lem:temporal-galerkin-solution} with $\beta=0$,  \eqref{eq:embedding-equatlity-I} and \eqref{eq:local-condition}, we infer
\begin{align}\label{eq:F(v(t))-F(v(s))}
\begin{split}
\|
P\big(
&F(X^{n}(t))-F(X^{n}(s))
\big)
\|_{L^p(\Omega;H)}
\\&
\leq
C
\|X^{n}(t)-X^{n}(s)\|_{L^{2p}(\Omega;H)}
\Big(
1+
\big(
\sup_{s\in[0,T]}
\|X^{n}(s)\|_{L^{2p}(\Omega;H^{\delta_0})}
\big)
\Big)^2
\\&
\leq
C|t-s|^{\frac{\delta_0}4},
\end{split}
\end{align}
and
\begin{align}\label{eq:bound-F(v(t))-H}
\begin{split}
\sup_{n\in \mathbb{N}}\sup_{s\in[0,T]}\|&PF(X^{n}(s))\|_{L^p(\Omega;H^1)}
\leq
C\sup_{n\in \mathbb{N}}\sup_{s\in[0,T]}\|F(X^{n}(s))\|_{L^p(\Omega;{\color{black}{H^1(D))}}}
\\
\leq
&
{\color{black} C }
\sup_{n\in \mathbb{N}}\sup_{s\in[0,T]}\big(\|\nabla F(X^{n}(s))\|_{L^p(\Omega;H)}+\|F(X^{n}(s))\|_{L^p(\Omega;H)}\big)
\\
\leq
&
C\sup_{n\in \mathbb{N}}\sup_{s\in[0,T]}\big(\big\|\|X^{n}(s)\|_1\|X^{n}(s)\|_V^2
\big\|_{L^{p}(\Omega;\mathbb{R})}
+\|X^{n}(s)\|_{L^p(\Omega;H^1)}+\|X^{n}(s)\|^3_{L^{3p}(\Omega;V)}\big)
\\
\leq
&
C\sup_{n\in \mathbb{N}}\sup_{s\in[0,T]}\big(\|X^{n}(s)\|_{L^{3p}(\Omega;H^1)}+
\|X^{n}(s)\|^3_{L^{3p}(\Omega;H^{\delta_0})}\big)
<\infty.
\end{split}
\end{align}
Then a combination of \eqref{I-spatio-temporal-S(t)}, \eqref{IV-spatio-temporal-S(t)}, \eqref{eq:F(v(t))-F(v(s))}
with \eqref{eq:bound-F(v(t))-H} shows for $\beta\in[0,\gamma]$
\begin{align}\label{eq:bound-integral-gamma}
\begin{split}
\Big\|&\int_s^tE(t-r)P_nAPF(X^{n}(r))\,\dd r\Big\|_{L^p(\Omega;H^\beta)}
\leq
\Big\|\int_s^tE(t-r)P_nAPF(X^{n}(t))\,\dd r\Big\|_{L^p(\Omega;H^\beta)}
\\
&
+
\int_s^t\big\|E(t-r)P_nAP(F(X^{n}(t))-F(X^{n}(r)))
\big\|_{L^p(\Omega;H^\beta)}\,\dd r
\\
\leq
&
C(t-s)^{\frac{3-\beta}4}\|PF(X^{n}(t))\|_{L^p(\Omega;H^1)}
+
\int_s^t(t-r)^{-\frac{2+\beta}4}
\big\|P\big(F(X^{n}(t))-F(X^{n}(r)))\big)
\big\|_{L^p(\Omega;H)}\,\dd r
\\
\leq
&
C(t-s)^{\frac{3-\beta}4}\sup_{n\in \mathbb{N}}\sup_{s\in[0,T]}\|PF(X^{n}(s))\|_{L^p(\Omega;H^1)}
+
\int_s^t(t-r)^{\frac{\delta_0-2-\beta}4}\,\dd r
\\
\leq
&
C(t-s)^{\frac{\gamma-\beta}4}.
\end{split}
\end{align}
{\color{black}{Eventually, the above estimate together with  \eqref{I-spatio-temporal-S(t)} and \eqref{II-spatio-temporal-S(t)}
leads to, for $\gamma\in(\delta_0,3]$ and any $\beta\in [0,\gamma]$}}
{\color{black}{\begin{align}
\begin{split}
\|v^{n}(t)\|_{L^p(\Omega;H^\gamma)}
&
\leq
\|E(t)P_nPX_0\|_{L^p(\Omega;H^\gamma)}
+
\Big\|
\int_0^t
E(t-s) P_n A P F(X^{n}(s))
\,\dd s
\Big\|_{L^p(\Omega;H^\gamma)}
\\
&
\leq
\|X_0\|_{L^p(\Omega;H^\gamma)}
+
C
\leq
C(T,\gamma,X_0,p)
\end{split}
\end{align}
and
\begin{align}
\begin{split}
\|v^{n}(t)-v^{n}(s)\|_{L^p(\Omega;H^\beta)}
\leq
&
\|(E(t-s)-I)v^{n}(s)\|_{L^p(\Omega;H^\beta)}
+
\Big\|
\int_s^t
E(t-r) P_n A P F(X^{n}(r))
\,\dd r
\Big\|_{L^p(\Omega;H^\beta)}
\\
\leq
&
C(t-s)^{\frac{\gamma-\beta}4} \|v^{n}(s)\|_{L^p(\Omega;H^\gamma)}
+
C(t-s)^{\frac{\gamma-\beta}4}.
\\
\leq
&
C(T,\gamma,X_0,p)(t-s)^{\frac{\gamma-\beta}4}.
\end{split}
\end{align}}}
Hence this finishes the proof of this theorem. $\square$

{\color{black}{At the moment, we are ready to deduce the strong convergence rate of the spectral Galerkin approximation.
\begin{theorem}\label{eq:strong-covergence-galerkin}
Let $v^n(t)$ be the solution of \eqref{eq:gerlrkin-approximation-CHC-i}
and
$X^n(t)=v^n(t)+ {\color{black}{Z^{n}(t)}}$.
If Assumptions \ref{assum:linear-operator-A}-\ref{assum:intial-value-data} are valid, then
\begin{align}\label{eq:theorem-strong-covergence-galerkin}
\sup_{s\in[0,T]}\|X^n(s)-X^m(s)\|_{L^p(\Omega;H)}
\leq
C \lambda_n^{-\frac\gamma2},\;m,n\in \mathbb{N}^+,\;m\geq n.
\end{align}
\end{theorem}}}
{\it Proof of Theorem \ref{eq:strong-covergence-galerkin}.}
 {\color{black}{From \eqref{eq:property-Pm} and \eqref{lemma:p-th-moment-bound-mild-solution}, it follows that, for any $\beta\in[0,\gamma)$ and $m\geq n$
\begin{align}
\begin{split}
\lim_{n\rightarrow\infty}\sup_{t\in[0,T]}\|Z^m(t)-Z^{n}(t)
\|_{L^p(\Omega;H^\beta)}
=
&
\lim_{n\rightarrow\infty}\sup_{t\in[0,T]}\|P^m(I-P_n) Z(t)\|_{L^p(\Omega;H^\beta)}
\\
\leq
&
C\lim_{n\rightarrow\infty} \lambda_{n+1}^{-\frac{\gamma-\beta}2}\sup_{t\in[0,T]}\|Z(t)\|_{L^p(\Omega;H^\gamma)}.
\end{split}
 \end{align}}}
Hence, it suffices to bound the error $\|v^m(t)-v^n(t)\|_{L^p(\Omega;H)}$.
 By  introducing an auxiliary problem
\begin{align}\label{eq:auxiliary-problem}
{\color{black}{\widetilde{v}^{n,m}(t)}}
=
E(t)P_nPX_0
-
\int_0^tE(t-s)P_nAPF(X^m(s))\,\dd s, n<m, n, m\in \mathbb{N},
\end{align}
we separate the error term  $\|v^m(t)-v^{n}(t)\|_{L^p(\Omega;H)}$ as follows
\begin{align}\label{eq:decomposement-error}
\|v^m(t)-v^{n}(t)\|_{L^p(\Omega;H)}
\leq
\|v^m(t)-\widetilde{v}^{n,m}(t)\|_{L^p(\Omega;H)}
+
\|\widetilde{v}^{n,m}(t)-v^{n}(t)\|_{L^p(\Omega;H)},
\quad n<m.
\end{align}
Resorting to \eqref{I-spatio-temporal-S(t)} {\color{black}{and the similar skills used in the proof of \eqref{eq:bound-integral-gamma}}}, we acquire that there exists a constant $C(X_0,\gamma,T)$ such that for  $t\in[0,T]$,
\begin{align}\label{eq:bound-auxiliary-problem}
\begin{split}
\sup_{m\geq n}\|\widetilde{v}^{n,m}(t)\|_{L^p(\Omega;H^\gamma)}
\leq
&
\|E(t)P_nPX_0\|_{L^p(\Omega;H^\gamma)}
+
\sup_{m\geq n}\Big\|\int_0^t
E(t-s)P_nAPF(X^m(s))\,\dd s\Big\|_{L^p(\Omega;H^\gamma)}
\\
\leq
&C \|X_0\|_{L^p(\Omega;H^\gamma)}
+
C
\leq
C(X_0,\gamma,T).
\end{split}
\end{align}
By letting ${\color{black}{\widetilde{X}^{n,m}(t)}}:=\widetilde{v}^{n,m}(t)+{\color{black}{Z^{n}(t)}}$,
the above estimate together with \eqref{lemma:p-th-moment-bound-mild-galerkin-solution}
shows
\begin{align}\label{eq:bound-mild-auxiliary-problem}
{\color{black}{\sup_{n\in \mathbb{N}}\sup_{m>n}}}
\sup_{s\in[0,T]}
\|
\widetilde{X}^{m,n}(s)
\|_{L^p(\Omega;H^\gamma)}
<\infty.
\end{align}
We are now ready to bound the first error item in \eqref{eq:decomposement-error}. In view of \eqref{eq:property-Pm}, \eqref{eq:bound-F-H} and
\eqref{I-spatio-temporal-S(t)}, we derive
\begin{align}\label{eq:convergence-rate-auxiliary-problem}
\begin{split}
\|v^m&(t)-\widetilde{v}^{n,m}(t)\|_{L^p(\Omega;H)}
\\&
\leq
\|E(t) P_m (I-P_n)P X_0\|_{L^p(\Omega;H)}
+\!
\|P_m (I-P_n)A^{-\frac\gamma2}\|_{\mathcal{L}(H)}
\Big\|\int_0^tE(t-s)A  P F(X^m(s))\,\dd s\|_{L^p(\Omega;H)}
\\&
\leq
C \lambda_{n+1}^{-\frac \gamma2}
\big(
\|X_0\|_{L^p(\Omega;H^{\gamma})}
+
 C
\big)
 \\&
\leq
C(X_0,T,\gamma) \lambda_{n+1}^{-\frac\gamma2}.
\end{split}
\end{align}
Next we turn our attention to the error ${\color{black}{\widetilde{e}^{n,m}(t)}}:=\widetilde{v}^{n,m}(t)-v^{n}(t)=\widetilde{X}^{n,m}(t)-X^{n}(t)$,
which is time differentiable and satisfies
\begin{align}\label{eq:error-equation-II}
\frac{\dd }{\dd t}\widetilde{e}^{n,m}(t)
+
 A^2 \widetilde{e}^{n,m}(t)
 =
 - P_n AP (F(X^m(t))-F(X^{n}(t))).
\end{align}
Taking inner product of \eqref{eq:error-equation-II} with $A ^{-1}\widetilde{e}^{n,m}(t)$ in $\dot{H}$  and making use of \eqref{eq:local-condition}, the fact
$
\|\widetilde{e}^{n,m}\|^2\leq |\widetilde{e}^{n,m}|_1|\widetilde{e}^{n,m}|_{-1}
$, {\color{black}{$\widetilde{e}^{n,m}\in \dot{H}$}}
and
$\widetilde{e}^{n}(0)=0$
imply
\begin{align}\label{eq:galerkin-approxiamtion-devative}
\begin{split}
\frac12\frac{\dd }{\dd s}|\widetilde{e}^{n,m}(s)|_{-1}^2
+
| \widetilde{e}^{n,m}(s)|_1^2
=
&
\big<
F(\widetilde{X}^{n}(s))-F(X^m(s))
,
\widetilde{e}^{n,m}(s)
\big>
+
\big<
F(X^{n}(s))-F(\widetilde{X}^{n,m}(s))
,
\widetilde{e}^{n,m}(s)
\big>
\\
\leq
&
\frac12
\Big\|
F(\widetilde{X}^{n,m}(s))-F(X^m(s))
\Big\|^2
+
\frac32
\|\widetilde{e}^{n,m}(s)\|^2
\\
\leq
&
C\|X^m(s)-\widetilde{X}^{n,m}(s)\|^2
\big(
1+
\|\widetilde{X}^{n,m}(s)\|_V^4
+
\|X^m(s)\|^4_V
\big)
\\
&
+
\frac12|\widetilde{e}^{n,m}(s)|_1^2
+
\frac98 |\widetilde{e}^{n,m}(s)|_{-1}^2
,
\end{split}
\end{align}
so that, after integration,
\begin{align}\label{eq:galerkin-approxiamtion-integrand}
|\widetilde{e}^{n}(t)|_{-1}^2
+
\int_0^t
|\widetilde{e}^{n}(s)|^2_1\,\dd s
\leq
C\int_0^t\|X^m(s)-\widetilde{X}^{n}(s)\|^2(1+\|\widetilde{X}^{n}(s)\|_V^4+
\|X^m(s)\|^4_V)\,\dd s,
\end{align}
where we also use  the Gronwall inequality.
By utilizing
\eqref{eq:bound-mild-auxiliary-problem},  \eqref{lemma:spatial-galerkin-solution}, \eqref{eq:convergence-rate-auxiliary-problem}
and \eqref{eq:embedding-equatlity-I} and applying the fact
\begin{align}
{\color{black}{\|Z^m(t)-Z^{n}(t)\|_{L^p(\Omega;H)}
\leq
\|P_m(I-P_n)Z(t)\|_{L^p(\Omega;H)}
\leq
C \lambda_{n+1}^{-\frac\gamma2}
\sup_{t\in[0,T]}\|Z(t)
\|_{L^p(\Omega;H^\gamma)},}}
\end{align}
one can find that,
\begin{align}\label{eq:axiliary-error-bound-H1-norm-integrand}
\begin{split}
\Big\|
\int_0^t
|\widetilde{e}^{n}(s)|_1^2
\,\dd s
\Big\|_{L^p(\Omega;\mathbb{R})}
\leq
&
C \Big\|
\int_0^t
\|X^m(s)-\widetilde{X}^{n}(s)\|^2
(1+
\|\widetilde{X}^{n}(s)\|_V^4
+
\|X^{m}(s)\|^4_V)
\,\dd s
\Big\|_{L^p(\Omega;\mathbb{R})}
\\
\leq
&
C(T)
\int_0^t
\Big(
\|v^m(s)-\widetilde{v}^{n}(s)\|_{L^{4p}(\Omega;H)}^2
+
{\color{black}{\|Z^m(t)-Z^{n}(t)\|_{L^{4p}(\Omega;H)}^2}}
\Big)
\,\dd s
\\
&\quad
\Big(
1+
\sup_{n\in \mathbb{N}}
\sup_{s\in[0,T]}
\|\widetilde{X}^{n}(s)\|_{L^{8p}(\Omega;V)}^4
+
\sup_{m\in \mathbb{N}}
\sup_{s\in[0,T]}
\|X^{m}(s)\|_{L^{8p}(\Omega;V)}^4
\Big)
\\
\leq
&
C(T,p,X_0,\gamma)\lambda_{n+1}^{-\gamma}.
\end{split}
\end{align}
At the moment we employ the above estimate to bound  $\|\widetilde{e}^{n}\|_{L^p(\Omega;H)}$, which can be  split  into two terms:
\begin{align}\label{eq:decompose-en}
\begin{split}
\|\widetilde{e}^{n}(t)\|_{L^p(\Omega;H)}
&\leq
\Big\|\int_0^tE(t-s)P_nAP\big[F(X^m(s))-F(X^{n}(s))\big]\,\dd s\Big\|_{L^p(\Omega;H)}
\\
&\leq
\Big\|\int_0^tE(t-s)P_nAP\big[F(X^m(s))-F(\widetilde{X}^{n}(s))\big]\,\dd s\Big\|_{L^p(\Omega;H)}
\\
&\quad +
\Big\|\int_0^tE(t-s)P_nAP\big[F(\widetilde{X}^{n}(s))-F(X^{n}(s))\big]\,\dd s\Big\|_{L^p(\Omega;H)}
\\
&=:
\text{Err}_1 + \text{Err}_2.
\end{split}
\end{align}
In the same way as in  \eqref{eq:axiliary-error-bound-H1-norm-integrand}, we obtain, $\kappa=\min\{\gamma,\frac d2+\frac14\}$
\begin{align}\label{eq:J_1-Galerkin-approximation}
\begin{split}
\text{Err}_1
&\leq
C
\int_0^t
(t-s)^{-\frac12}
\|F(X^m(s))-F(\widetilde{X}^{n}(s))\|_{L^p(\Omega;H)}
\,\dd s
\\
&\leq
C
\int_0^t
(t-s)^{-\frac12}
\|X^m(s)-\widetilde{X}^{n}(s)\|_{L^{2p}(\Omega;H)}
\,\dd s
\\
&
\qquad\quad
\sup_{m,n\in \mathbb{N},s\in[0,T]}
\Big(
1+\|X^m(s)\|_{L^{4p}(\Omega;V)}^2
+\|\widetilde{X}^{n}(s)\|_{L^{4p}(\Omega;V)}^2
\Big)
\\
&
\leq
C(T,p,X_0,\gamma)
\lambda_{n+1}^{-\frac\gamma2}.
\end{split}
\end{align}
To bound the term $\mathrm{Err}_2$, we apply \eqref{eq:FX1-FX2-d=1}  \eqref{eq:FX1-FX2-d=23}  \eqref{eq:bound-mild-auxiliary-problem}, \eqref{eq:axiliary-error-bound-H1-norm-integrand}, \eqref{lemma:spatial-galerkin-solution} and  the fact $\widetilde{e}^{n}\in \dot{H}^1$  to derive that, for $\eta=\min\{\gamma,\frac34\}$ and  $\kappa=\min\{\gamma, \frac {2d+1}4\}$
\begin{align}
\begin{split}
\text{Err}_2
\leq
&
C
\Big\|
\int_0^t
(t-s)^{-\frac{2-\eta}4}
\|F(\widetilde{X}^{n}(s))-F(X^{n}(s))\|_{\eta}
\,\dd s
\Big\|_{L^p(\Omega;\mathbb{R})}
\\
\leq
&C
\Big\|
\int_0^t
(t-s)^{-\frac{2-\eta}4}
|\widetilde{e}^{n}(s)|_1
\big(
1
+
\|\widetilde{X}^{n}(s)\|_\kappa^2
+
\|X^{n}(s)\|_\kappa^2
\big)
\,\dd s
\Big\|_{L^p(\Omega;\mathbb{R})}
\\
\leq
&
C
\Big\|
\int_0^t
|\widetilde{e}^{n}(s)|_1^2
\,\dd s
\Big\|_{L^p(\Omega;\mathbb{R})}^{\frac12}
\Big\|
\int_0^t
(t-s)^{-\frac{2-\eta}2}
\big(
1+\|\widetilde{X}^{n}(s)\|_\gamma^4+\|X^{n}(s)\|_\gamma^4
\big)\,\dd s
\Big\|^{\frac12}_{L^p(\Omega;\mathbb{R})}
\\
\leq
&
C(T,X_0,\gamma,p)
\lambda_{n+1}^{-\frac\gamma2}.
\end{split}
\end{align}
Finally, gathering the estimate of $E_1$ and $E_2$ together gives
\begin{align}
\|v^{n}(t)-\widetilde{v}^{n,m}(t)\|_{L^p(\Omega;H)}
\leq
C\lambda_{n+1}^{-\frac\gamma2},
\end{align}
which in combination with  \eqref{eq:convergence-rate-auxiliary-problem} shows
\begin{align}\label{eq:v-v-in-h}
{\color{black}{\sup_{m\geq n}\sup_{t\in[0,T]}}}\|v^m(t)-v^{n}(t)\|_{L^p(\Omega;H)}
\leq
C\lambda_{n+1}^{-\frac\gamma2}.
\end{align}
Hence this  finishes the proof. $\square$
\subsection{Existence, uniqueness and regularity of the mild solution}
This part is devoted  to the well-posedness and spatio-temporal  regularities of the mild solution to the problem \eqref{eq:abstract-CHC-form}.
\begin{theorem}\label{them:regulairty-mild-solution}
If Assumptions \ref{assum:linear-operator-A}-\ref{assum:intial-value-data} are valid, then the problem \eqref{eq:abstract-CHC-form} admits a unique mild solution given by
\begin{align}\label{them:eq-mild-solution-stochastic-equation}
X(t)
=
E(t)X_0
-
\int_0^t
E(t-s) A P F(X(s))
\,\mathrm{d} s
+
\int_0^t
E(t-s)
\mathrm{d} W(s).
\end{align}
Moreover, for 
any $p\geq 1$ we have
\begin{align}\label{them:spatial-regularity-mild-stoch}
\sup_{s\in[0,T]}\|X(s)\|_{L^p(\Omega;H^\gamma)}
<\infty,
\end{align}
and for any $\beta\in[0,\gamma]$, there exists a constant $C$ depending on $T$, $p$ such that
\begin{align}\label{them:temporal-regularity-mild-stoch}
\|X(t)-X(s)\|_{L^p(\Omega;H^\beta)}
\leq
C(t-s)^{\min\{\frac12,\frac{\gamma-\beta}4\}}.
\end{align}
\end{theorem}
  {\it Proof of Theorem \ref{them:regulairty-mild-solution}.} By the same arguments used in the proof of Lemma \ref{eq:approximate-solution-higher-norm}, we can show \eqref{them:spatial-regularity-mild-stoch} and \eqref{them:temporal-regularity-mild-stoch}.
Hence we only  prove   the existence  and  uniqueness of the mild solution to  \eqref{eq:abstract-CHC-form}.
{\color{black}{From \eqref{eq:v-v-in-h}, \eqref{lemma:spatial-galerkin-solution} and  the Sobolev interpolation inequality, if follows that
\begin{align}
\begin{split}\label{eq:v-v-in-h14}
\| X^m(t)-X^{n}(t)\|_{L^p(\Omega;H^{\frac14})}
\leq&
C\|X^m(t)-X^{n}(t)\|^{\frac12}_{L^{2p}(\Omega;H)}\|X^m(t)-X^{n}(t)\|^{\frac12}_{L^{2p}(\Omega;H^{\frac12})}
\\
\leq&
C\|X^m(t)-X^{n}(t)\|^{\frac12}_{L^{2p}(\Omega;H)}\Big(\sup_{n\geq 0}\sup_{t\in[0,T]}\|X^n(t)\|_{L^{2p}(\Omega;H^\gamma)}\Big)^{\frac12}
\\
\leq
&
C\lambda_{n+1}^{-\frac\gamma 4}\rightarrow 0.
\end{split}
\end{align}
 By proceeding with the same process of the above proof and   using bootstrapping arguments, we can prove
\begin{align}
\sup_{m\geq n}\sup_{t\in[0,T]}\|X^m(t)-X^{n}(t)\|_{L^p(\Omega;H^{\frac d3})}\rightarrow0
\end{align}
This shows that the sequence ${\color{black}{\{X^{n}(\cdot)\}}}$ is Cauchy in $C_W([0,T];H^{\frac d3}), d=1,2,3$.
Let
\begin{align}
X ( \cdot )
=
\lim_{n\rightarrow\infty}
X^{n}( \cdot )
\;
\text{ in }
\;
C_W([0,T];H^{\frac d3}), d=1,2,3.
\end{align}
Then, the following results hold, for $d=1,2,3$
\begin{align}\label{eq:bound-X}
\sup_{t\in[0,T]}\|X(t)\|_{L^p(\Omega;H^{\frac d3})}<\infty,
\end{align}
 and
\begin{align}\label{eq:convergence-Xn-X}
\lim_{n\rightarrow\infty}
X^{n}( \cdot ) = X ( \cdot )
\;
\text{ in }
\;
C_W([0,T];H^\frac d3).
\end{align}
In what follows, we will prove that $X(t)$ is a mild solutions of \eqref{eq:abstract-CHC-form-I}.
It suffices to show that
\begin{align}\label{eq:conergence-L}
\lim_{n\rightarrow \infty}\underbrace{\left\|\int_0^t
E(t-s) P_n A P
F(
X^{n}(s)
)
\,\dd s
-
\int_0^t
E(t-s)  A P
F(
X(s)
)
\,\dd s
\right\|_{L^p(\Omega;H)}}_{\mathbb{L}}\rightarrow0.
\end{align}
To show the above result,  we split $\mathbb{L}$ into the following two parts
\begin{align}
\begin{split}\label{eq:L-decompose}
\mathbb{L}
\leq
&
\Big\|\int_0^t
E(t-s)  A(P_n -I)P
F(
X^{n}(s)
)
\,\dd s\Big\|_{L^p(\Omega;H)}
\\
&
+
\Big\|\int_0^t
E(t-s)  A P\big(F(
X^{n}(s)
)-
F(
X(s)
)\big)
\,\dd s
\Big\|_{L^p(\Omega;H)}
=\mathbb{L}_1+\mathbb{L}_2.
\end{split}
\end{align}
{\color{black}{By using \eqref{eq:approximate-solution-higher-norm}, \eqref{I-spatio-temporal-S(t)} with $\mu=
\frac34$ and \eqref{eq:property-Pm} with $\alpha=\frac12$,  the term $\mathbb{L}_1$ can be estimated as follows
\begin{align}\label{eq:convergence-L1}
\begin{split}
\mathbb{L}_1
\leq
&
\int_0^t
\|E(t-s)  A^{\frac32}\|_{\mathcal{L}(H)}\|(I-P_n)A^{-\frac12}\|_{\mathcal{L}(H)}\|P
F(
X^{n}(s)
)\|_{L^p(\Omega;H)}
\,\dd s
\\
\leq
&
C\lambda_n^{-\frac12}\int_0^t(t-s)^{-\frac34}\,\dd s\big(1+\sup_{t\in[0,T]}\|X^n(t)\|_{L^{3p}(\Omega;L^6)}^3\big)
\\
\leq
&
C\lambda_n^{-\frac12}\rightarrow 0.
\end{split}
\end{align}}}
For the term $\mathbb{L}_2$, we use \eqref{eq:bound-X}, \eqref{lemma:spatial-galerkin-solution}, \eqref{eq:convergence-Xn-X}, \eqref{eq:embedding-equatlity-III}, \eqref{eq:embedding-equatlity-I} and  \eqref{I-spatio-temporal-S(t)} with $\mu=
\frac34$ to infer
\begin{align}
\begin{split}
\mathbb{L}_2
\leq
&
\int_0^t
\big\|
E(t-s)  A^{\frac32}\|_{\mathcal{L}(H)}\|A^{-\frac12} P\big(F(
X^{n}(s)
)-
F(
X(s)
)\big)
\big\|_{L^p(\Omega;H)}
\,\dd s
\\
\leq
&C\int_0^t
\big\|
(t-s)^{-\frac34}\|F(
X^n(s)
)-
F(
X(s)
)
\big\|_{L^p(\Omega;L^{\frac65})}
\,\dd s
\\
\leq
&\int_0^t
(t-s)^{-\frac34}\Big\|\|
X^{n}(s)
-
X(s)
\|
\big(1+\|X^n(t)\|_{L^6}^2+\|X(t)\|_{L^6}^2\big)\Big\|_{L^p(\Omega;\mathbb{R})}
\,\dd s
\\
\leq
&C\int_0^t
(t-s)^{-\frac34}\|
X^{n}(s)
-
X(s)
\|_{L^{2p}(\Omega;H)}\dd s \\
\quad
&
\big(1+\sup_{t\in[0,T]}(\|X^n(t)\|_{L^{4p}(\Omega;L^6)}+\|X(t)\|_{L^{4p}(\Omega;L^6)})^2\big)
\\
\leq
&
C\sup_{s\in[0,T]}\|
X^{n}(s)
-
X(s)
\|_{L^{2p}(\Omega;H)}\rightarrow 0.
\end{split}
\end{align}
This together with \eqref{eq:convergence-L1} and \eqref{eq:L-decompose} shows \eqref{eq:conergence-L}.}}
The existence of the mild solution is proved.

Next, let us prove the uniqueness. Let $X_1$ and $X_2$ be two  mild solutions of \eqref{eq:abstract-CHC-form}. Then $X_1-X_2$ is the mild solution of the deterministic problem
\begin{align}
\left\{
\begin{array}{ll}
\frac{\mathrm{d}}{\mathrm{d}t}(X_1-X_2)
+ A
\big(
A(X_1-X_2)
+F(X_1)-F(X_2)
\big)=0,&
\\
(X_1-X_2)(0)
=0.&
\end{array}
\right.
\end{align}
Multiplying both sides of the first identity by $A^{-1}(X_1-X_2)$ and applying {\color{black}{the fact $X_1-X_2\in \dot{H}$}} and  the similar arguments as in \eqref{eq:galerkin-approxiamtion-devative}  yield
\begin{align}
\tfrac12
\tfrac{\mathrm{d}}{\mathrm{d}t}
|X_1-X_2|_{-1}^2
+
|X_1-X_2|_1^2
\leq
\tfrac12
|X_1-X_2|_1^2
+
\tfrac12
|X_1-X_2|_{-1}^2,
\end{align}
which implies $X_1=X_2$.
So uniqueness is shown and the proof of this theorem is complete. $\square$

\begin{lemma}\label{lem:spatial-regulairty-F-mild}
If Assumptions \ref{assum:linear-operator-A}-\ref{assum:intial-value-data} are valid,
the following results hold,
\begin{align}\label{lem:eq-F-mild-bound}
\sup_{s\in[0,T]}\|F(X(s))\|_{L^p(\Omega;H^{\ell(\gamma)})}<\infty,
\end{align}
and for $t \geq s$,
\begin{align}\label{lem:temporal-regulairty-F}
\|F(X(t))-F(X(s))\|_{L^p(\Omega;H^{l(\gamma)-2})}
\leq
 C(t-s)^{\frac18},
\end{align}
where
\begin{equation}
\label{eq:l-gamma-defn}
l(\gamma)=\left\{
\begin{array}{ll}
2,&for\;\gamma\in (\frac d 2,2],\\
3,&for\; \gamma\in(2,3],\\
4,&for\;\gamma\in (3,4],
\end{array}
\right.
\;and
\quad
\ell(\gamma)=\left\{
\begin{array}{ll}
0, &for\;\gamma\in(\frac d2,2],
\\
2,&for\;\gamma\in(2,4].
\end{array}
\right.
\end{equation}
\end{lemma}
 {\it Proof of Lemma \ref{lem:spatial-regulairty-F-mild}.}
 The result \eqref{lem:eq-F-mild-bound} can be easily shown by using  \eqref{eq:embedding-equatlity-I}, \eqref{eq:algebra-properties-Hs} and \eqref{them:spatial-regularity-mild-stoch}.
  Next we show \eqref{lem:temporal-regulairty-F}. The case $\gamma\in (\frac d2,2] $ can be proved  by utilizing \eqref{them:temporal-regularity-mild-stoch}, \eqref{eq:local-condition} and \eqref{eq:embedding-equatlity-I}. For the case $\gamma\in (2,3]$, we first utilize the same arguments used in  the proof of \eqref{eq:FX1-FX2-d=23} to obtain
 \begin{align}
 \|F(X(t))-F(X(s))\|_1
 \leq
 (1+\|X(t)\|_\gamma^2+\|X(s)\|_\gamma^2)\|X(t)-X(s)\|_1.
 \end{align}
This together with \eqref{them:spatial-regularity-mild-stoch} and   \eqref{them:temporal-regularity-mild-stoch} shows the case $\gamma\in(2,3]$. Similarly, by employing  \eqref{eq:algebra-properties-Hs} and \eqref{them:temporal-regularity-mild-stoch}, one can show the case $\gamma\in (3,4]$. Hence, the proof of this lemma is complete. $\square$

\section{The fully discrete approximation}
In this section, we consider a full discretization of the CHC equation and show the maximal-type  moment bounds of the solution to the fully discrete problem, which will be used later to the convergence analysis. Throughout the proofs, $C$ denotes a generic nonnegative constant that is independent of the discretization parameters $N$ and $k$ and may change from lines to lines.

 Let
 $k=T/M$, $M\in \mathbb{N}^+$
 be a time step-size and
 $t_m=mk,\, m\in\{0,1,\cdots,M\}$.
 For
 $N \in \mathbb{N}$,
 we define a finite dimensional subspace of $H$ by
 $H_N:= span\{{\color{black}{e_0}}, e_1, e_2, \cdots, e_N \}$
 and the projection
 $P_N:H^\alpha\rightarrow H_N$
 such that
 $P_N\xi=\sum_{j=0}^{N}\big<\xi,e_j\big>e_j$,
 for
 $\xi\in H^\alpha$,
 $\alpha\in \mathbb{R}$.
It is not difficult to deduce that
\begin{align}\label{eq:interpolation-property}
\|(I-P_N)\phi\|
\leq
C \lambda_N^{-\frac\alpha2}|\phi|_\alpha,\;\forall \phi\in H^\alpha,\; \alpha>0.
\end{align}
Then the fully discrete approximation of  the problem \eqref{eq:abstract-CHC-form} is to find $X_m^{M,N}\in H_N$ such that
\begin{align}\label{eq:full-discretization}
X_m^{M,N}-X_{m-1}^{M,N}
+
k A^2 X_m^{M,N}
+
kP_NAF(X_m^{M,N})
=
P_N \Delta W_m,\; X_0^{M,N}=P_NX_0,\;
\end{align}
where we write $\Delta W_m=W(t_m)-W(t_{m-1})$ for brevity.
By introducing a family of operators $\{{\color{black}{E_{k,N}^m}}\}_{m=1}^M$:
\begin{align}\label{eq:full-discrete-solution-operator}
E_{k,N}^m v=(I+kA^2)^{-m} P_N v
=
{\color{black}{\sum_{j=0}^{N}(1+k \lambda_j^2)^{-m}\big<v,e_j\big>e_j,\;\forall v\in H}},
\end{align}
the solution of \eqref{eq:full-discretization} can be written as
\begin{align}
X_m^{M,N}=E_{k,N}^mP_N X_0
-
k\sum_{j=1}^mE_{k,N}^{m-j+1}P_NAF(X_j^{M,N})
+
\sum_{j=1}^mE_{k,N}^{m-j+1}P_N\Delta W_j.
\end{align}
{\color{black}{Noting that the above implicit scheme works on the space $H_N$ and that the mapping
$A^2
+kP_NAF(\cdot)$ obeys a kind of monotonicity condition in the Hilbert space $\big(H_N, (\cdot, \cdot)_{-1}\big)$, 
one can see that the implicit scheme \eqref{eq:full-discretization} is well-posed in $\big(H_N, (\cdot, \cdot)_{-1}\big)$. }}

The next theorem offers a priori moment bounds for the fully discrete approximation.
\begin{theorem}\label{them:bound-numerical-solution} {\color{black}{Let $k=T/M, M\in \mathbb{N}^+$ be a time step-size and $X_m^{M,N}$ be the solution of \eqref{eq:full-discretization}.}}
 Suppose Assumptions
 \ref{assum:linear-operator-A}-\ref{assum:intial-value-data} are valid,
 then there
 exists a positive constants $k_0$ such that for all
 $k\leq k_0$, $N\in \mathbb{N}$
 and $\forall\, p\geq 1$,
\begin{align}\label{lem:eq-bound-solution-full-stochatic}
\sup_{1\leq m\leq M}
\|X_m^{M,N}\|_{L^p(\Omega;H^\kappa)}
<
\infty,
\;\kappa=\min\big\{\gamma,\tfrac d2+\tfrac 14\big\}.
\end{align}
\end{theorem}

Before showing it,   we need to introduce  some {\color{black}{smoothing}} properties of $E_{k£¬N}^m$,
 which {\color{black}{is a variant of \eqref{I-spatio-temporal-S(t)} and \eqref{III-spatio-temporal-S(t)} and}}  can be proved by a slight modification of the proof of \cite[Lemma 3.2]{ruisheng2018error}.

\begin{lemma}\label{lem:bound-rational-approximation-semi}
Under Assumption \ref{assum:linear-operator-A}, the following estimates for  $E_{k}^{m}$ hold.\\
(i) Let $\mu\in[0,2]$. There exists a constant $C$ such that
\begin{align}\label{lem:bounded-time-full-deterministic-problem}
\|A^\mu E_{k,N}^mv\|
&\leq
Ct_m^{-\frac\mu2}\|v\|,\quad\forall v\in H,
\;m= 1, 2,3,  \cdots, M.
\end{align}
(ii) There exists a constant $C$ such that
\begin{align}\label{lem:eq-bound-sum-full-operator-II}
\Big(k\sum_{j=1}^m\|A E_{k,N}^j v\|^2\Big)^{\frac12}
\leq
C\|v\|,\quad \forall v\in H.
\end{align}
\end{lemma}

{\color{black}{
In addition, we need the following lemma concerning a strong moment bounds of the discrete stochastic convolution.
To arrive at it, we extend  the classical factorization method used in the continuous setting \cite[Chapter 5]{da2014stochastic}.
\begin{lemma}
\label{lem:stoch-conv-strong-moment}
For any $\mu \in\left(\frac d2,\min\{\gamma,2\}\right)$, it holds that
\begin{align}
\mathbb{E}\Big[ \sup_{ 1\leq m\leq M } \Big \|  \sum_{i=1}^m E_{k,N}^{m-i+1}P_N \Delta W_i \Big\|_\mu^p \Big]
\leq
C ( \|A^{\frac{\gamma-2}2}PQ^{\frac12}\|^p_{\mathcal{L}_2(H)}+\|(I-P)Q^{\frac12}\|^p_{\mathcal{L}_2(H)})
<
\infty.
\end{align}
\end{lemma}
%
{\it Proof of Lemma \ref{lem:stoch-conv-strong-moment}.}
For $i \in \{ 1, 2, \cdots, M\}$,  we define 
\begin{equation}
\tilde{E}_{k,N}(t) = E_{k,N}^{i},
\quad
t \in (t_{i-1}, t_{i}].
\end{equation}
In this case $\tilde{E}_{k,N}(t)$ is not continuous with respect to $t$ and does not have the semigroup property.
However, we have a kind of weak semigroup property as follows:
\begin{equation}
\label{eq:weak-semigroup-prop}
\tilde{E}_{k,N}(t) 
= 
E_{k,N}^{i}
=
E_{k,N}^{j} E_{k,N}^{ i - j}
=
E_{k,N}^j \tilde{E}_{k,N}(t - t_j ),
\quad
0\leq t_j \leq t_{i-1} < t \leq t_{i}.
\end{equation} 
Also, we define a continuous version of the discrete stochastic convolution as
\begin{align}
\begin{split}
\mathcal{O}_{k,N}(t)
: =
\int_0^t \tilde{E}_{k,N}(t-s)\textmd{d}W(s).
\end{split}
\end{align}
It is not difficult to see that $\mathcal{O}_{k,N}(t_m) = \sum_{i=1}^m E_{k,N}^{m-i+1}P_N \Delta W_i $.
We shall first extend  the classical factorization method in the continuous setting \cite[Chapter 5]{da2014stochastic}, 
which is based on the following elementary identity
\begin{align}
\int_\sigma^t(t-s)^{\alpha-1}(s-\sigma)^{-\alpha}\,\dd s
=
\frac{\pi}{\sin \pi\alpha},\; 0 \leq \sigma \leq t,\;0<\alpha<1.
\end{align}
In what follows we suppose $\alpha\in(0,\frac12)$ and  $p_0>\frac1{2\alpha}$.
Using the identity and the stochastic Fubini theorem we obtain
\begin{align}
\mathcal{O}_{k,N}(t_m)
 = &
\frac{\sin \pi\alpha}{\pi}
\int_0^{t_m}\int_\sigma^{t_m} 
\tilde{E}_{k,N}
(t_m-\sigma)(t_m - s)^{\alpha-1}(s-\sigma)^{-\alpha}\,\dd s\,\dd W(\sigma)
\nonumber
\\
 = &
\frac{\sin \pi\alpha}{\pi}
\int_0^{t_m}\int_0^s
\tilde{E}_{k,N}(t_m-\sigma)( t_m - s)^{\alpha-1}(s-\sigma)^{-\alpha}\,\dd W(\sigma)\,\dd s
.
\end{align}
Further, one can recall \eqref{eq:weak-semigroup-prop} and deduce that
\begin{align}
\begin{split}
\mathcal{O}_{k,N}(t_m)
=
&
\frac{\sin \pi\alpha}{\pi}
\sum_{j=1}^m\int_{t_{j-1}}^{t_j}\int_0^s
\tilde{E}_{k,N}(t_m-\sigma)( t_m - s)^{\alpha-1}(s-\sigma)^{-\alpha}\,\dd W(\sigma)\,\dd s
\\
=
&
\frac{\sin \pi\alpha}{\pi}
\sum_{j=1}^m\int_{t_{j-1}}^{t_j}
E_{k,N}^{m-j }( t_m - s)^{\alpha-1}\int_0^s
\tilde{E}_{k,N}( t_{j} - \sigma) ( s - \sigma)^{-\alpha}\,\dd W(\sigma)\,\dd s
\\
=
&
\frac{\sin \pi\alpha}{\pi}
\sum_{j=1}^m
\int_{t_{j-1}}^{t_j}E_{k,N}^{m-j}(t_m - s)^{\alpha-1} \Theta_{\alpha} ( s )
\,\dd s,
\end{split}
\end{align}
where we denote $ \lceil  s \rceil = t_{j}$ for $s \in (t_{j-1}, t_j]$ and 
\begin{align}
\Theta_{\alpha}(s)
: =
\int_0^s
\tilde{E}_{k,N}( \lceil  s \rceil - \sigma)(s-\sigma)^{-\alpha}\,\dd W(\sigma).
\end{align}
By the H\"{o}lder inequality and the stability of $E_{k, N}$, we derive that,  for $p_0>\frac1{2\alpha}$ with $\alpha\in(0,\frac12)$,
\begin{align}
\begin{split}
\|\mathcal{O}_{k,N}(t_m)\|_\mu^{2p_0}
\leq
&
C\left(\sum_{j=1}^m\int_{t_{j-1}}^{t_j}\|E_{k,N}^{m-j}
\|^{\frac{2p_0}{2p_0-1}}_{\mathcal{L}(H)}( t _m- s)^{\frac{2p_0(\alpha-1)}{2p_0-1}}\,\dd s
\right)^{2p_0-1}\int_0^{t_m}\| \Theta_{\alpha} ( s ) \|_\mu^{2p_0}\,\dd s
\\
\leq
&
C\int_0^{t_m}\| \Theta_{\alpha} ( s ) \|_\mu^{2p_0}\,\dd s.
\end{split}
\end{align}
Using the Burkholder-Davis-Gundy-type inequality and \eqref{lem:bounded-time-full-deterministic-problem}, 
we derive that, for $\mu \in\left(\frac d2,\min\{\gamma,2\}\right)$, 
\begin{align}
\begin{split}
\mathbb{E}
& \Big[
\sup_{1\leq m\leq M}\| \mathcal{O}_{k,N}(t_m)\|_\mu^{2p_0}
\Big]
\leq
C\int_0^T\mathbb{E}\left[\| \Theta_{\alpha} ( s ) \|_\mu^{2p_0}\right]\,\dd s
\\
& \quad
\leq
C\int_0^T\left(\int_0^s
(s-\sigma)^{-2\alpha}
\Big(
\|A^{\frac\mu2}
\tilde{E}_{k,N}( \lceil  s \rceil - \sigma)Q^{\frac12}\|_{\mathcal{L}_2(H)}^2
+
\|(I-P)Q^{\frac12}\|^{2}_{\mathcal{L}_2(H)}
\Big)\,\dd \sigma\right)^{p_0}\,\dd s
\\
& \quad
\leq
C\int_0^T\left(\int_0^s
(s-\sigma)^{-2\alpha}
\Big(
(s-\sigma)^{-\frac{\max\{\mu-\gamma+2,0\}}2}
\|A^{\frac{\gamma-2}2} P Q^{\frac12}\|_{\mathcal{L}_2(H)}^{2}
+
\|(I-P)Q^{\frac12}\|^{2}_{\mathcal{L}_2(H)}
\Big)\,\dd \sigma\right)^{p_0}\,\dd s
\\
& \quad
\leq
C\int_0^T\left(\int_0^s
(s-\sigma)^{-2\alpha-\frac{\max\{\mu-\gamma+2,0\}}2}\,\dd \sigma \right)^{p_0}
\Big(
\|A^{\frac{\gamma-2}2} P Q^{\frac12}\|_{\mathcal{L}_2(H)}^{2p_0}
+
\|(I-P)Q^{\frac12}\|^{2p_0}_{\mathcal{L}_2(H)}
\Big)\,
\dd s
<\infty,
\end{split}
\end{align}
where we require that $\alpha<\frac12-\frac{\max\{\mu-\gamma+2,0\}}4$. $\square$
}}

At the moment, we are ready to prove Theorem \ref{them:bound-numerical-solution}.

{\it Proof of Theorem \ref{them:bound-numerical-solution}.}
{\color{black}{We begin with showing the momet boundedness of}} 
\begin{equation*}
Z_m^{M,N}
:=
E_{k,N}^m{\color{black}{P_N}}(I-P)X_0
+
\sum_{i=1}^m E_{k,N}^{m-i+1}P_N \Delta W_i.
\end{equation*}
{\color{black}{
%
%
%
Thanks to Lemma \ref{lem:stoch-conv-strong-moment} and applying \eqref{lem:bounded-time-full-deterministic-problem} with $\mu=0$ give
\begin{align}\label{eq:bound-full-discrete-stochastic-convulution}
\begin{split}
\mathbb{E} \Big[ \sup_{ 1\leq m\leq M } \|Z_m^{M,N} \|_\gamma^p\Big]
\leq
&
C\mathbb{E}\Big[ \sup_{ 1\leq m\leq M } \|E_{k,N}^m{\color{black}{P_N}}(I-P)X_0\|_\gamma^p\Big]
+
\mathbb{E}\Big[ \sup_{ 1\leq m\leq M } \Big \|  \sum_{i=1}^m E_{k,N}^{m-i+1}P_N \Delta W_i \Big\|_\gamma^p \Big]
\\
\leq
&
C\mathbb{E}\big[\|X_0\|_\gamma^p\big]
+
C ( \|A^{\frac{\gamma-2}2}PQ^{\frac12}\|^p_{\mathcal{L}_2(H)}+\|(I-P)Q^{\frac12}\|^p_{\mathcal{L}_2(H)})
<\infty,
\end{split}
\end{align}
}}
Similarly to the continuous case, we only need to consider the spatial regularity of the solution to the following problem
\begin{align}
Y_m^{M,N}:
=
X_m^{M,N}-{\color{black}{Z_m^{M,N}}}
=
E_{k,N}^mPX_0
-
k\sum_{i=1}^{m}E_{k,N}^{m+1-i}AF(Y_i^{M,N}+{\color{black}{Z_i^{M,N}}}).
\end{align}
It is straightforward to verify that $Y_m^{M,N}$ for $m\in\{1,2,\cdots,M\}$ satisfies
\begin{align}\label{eq:full-modified-problem}
\frac{Y_m^{M,N}-Y_{m-1}^{M,N}}k
+
A^2 Y_m^{M,N}
=
-P_NAP F(Y_m^{M,N}+{\color{black}{Z_m^{M,N}}}),\;Y_0^{M,N}=P_N{\color{black}{P}} X_0.
\end{align}
Multiplying both sides of \eqref{eq:full-modified-problem} by $A^{-1}Y_m^{M,N}$ and applying the same arguments used in the proof of \eqref{eq:vm-H-1-norm},  we have
\begin{align}
\begin{split}
|Y_m^{M,N}|_{-1}^2
\leq
&
|Y_{m-1}^{M,N}|_{-1}^2
-
2k|Y_m^{M,N}|_1^2
-
2k\big<F(Y_m^{M,N}+{\color{black}{Z_m^{M,N})}}, Y_m^{M,N}
\big>
\\
=
&
|Y_{m-1}^{M,N}|_{-1}^2
-
2k|Y_m^{M,N}|_1^2
-
2k\big<(Y_m^{M,N})^3-Y_m^{M,N}, Y_m^{M,N}
\big>
\\
&
-
2k\big<{\color{black}{3(Y_m^{M,N})^2 Z_m^{M,N}+3Y_m^{M,N} (Z_m^{M,N})^2+(Z_m^{M,N})^3-Z_m^{M,N}}},Y_m^{M,N}\big>
\\
\leq
&
|Y_{m-1}^{M,N}|_{-1}^2
-
2k|Y_m^{M,N}|_1^2
-
2k\|Y_m^{M,N}\|_{L_4}^4
+
2k \|Y_m^{M,N}\|^2
\\
&
+
2k \Big( 3\|Y_m^{M,N}\|_{L^4}^3\|Z_m^{M,N}\|_{L^4}
+
3\|Y_m^{M,N}\|^2_{L^4}\|Z_m^{M,N}\|^2_{L^4}
\\ 
&
+
\|Y_m^{M,N}\|_{L^4}\|Z_m^{M,N}\|^3_{L^4}
+
\|Y_m^{M,N}\|_{L^2}\|Z_m^{M,N}\|_{L^2}
\Big)
\\
\leq
&
|Y_{m-1}^{M,N}|_{-1}^2
-
2k|Y_m^{M,N}|_1^2
-
k\|Y_m^{M,N}\|_{L^4}^4
+
Ck(1+\|{\color{black}{Z_m^{M,N}}}\|_{L^4}^4).
\end{split}
\end{align}
After repeated application, this yields
\begin{align}\label{eq:bound-L4-summation}
|Y_m^{M,N}|_{-1}^2
+
2\sum_{j=1}^mk|Y_j^{M,N}|_1^2
+
\sum_{j=1}^mk\|Y_j^{M,N}\|_{L^4}^4
\leq
C+\|X_0\|_{-1}^2
+
Ck\sum_{j=1}^m\|{\color{black}{Z_j^{M,N}}}\|_{L^4}^4,
\end{align}
which together with \eqref{eq:bound-full-discrete-stochastic-convulution}   yields
\begin{align}\label{eq:bound-L4-summation-expectation}
\begin{split}
\mathbb{E}\Big[\Big(k\sum_{j=1}^M|Y_j^{M,N}|_1^2\Big)^p\Big]
+
\mathbb{E}\Big[\Big(k\sum_{j=1}^M\|Y_j^{M,N}\|_{L^4}^4\Big)^p\Big]
\leq
&
C(T,p, X_0)\big(1+\mathbb{E}\Big[\Big(k\sum_{j=1}^M\|{\color{black}{Z_j^{M,N}}}\|_{L^4}^4\Big)^p\Big]\big)
\\
\leq
&
C(T,p, X_0)\Big(1+k\sum_{j=1}^M\mathbb{E}
\big[\|{\color{black}{Z_j^{M,N}}}\|_\kappa^{4p}
\big]\Big)
<\infty.
\end{split}
\end{align}
Next we consider the bound of $\|Y_m^{M,N}\|$.
By taking the inner product on both sides of \eqref{eq:full-modified-problem}
with $Y_{m}^{M,N}$ and using integration by parts formula,
we have
\begin{equation}
\begin{split}
\|Y_{m}^{M,N}\|^2
\leq
&
\|Y_{m-1}^{M,N}\|^2
-
2k
\|AY_{m}^{M,N}\|^2
-
2k\big<P_NAF(Y_{m}^{M,N}+{\color{black}{Z_{m}^{M,N}}}), Y_{m}^{M,N}\big>
\\
\leq
&
\|Y_{m-1}^{M,N}\|^2
-
2k\|AY_{m}^{M,N}\|^2
-
\tfrac32 k \|\nabla( Y_{m}^{M,N})^2\|^2
+
2k \|\nabla Y_{m}^{M,N}\|^2
\\
&
-
2k\big<{\color{black}{3(Y_{m}^{M,N})^2 Z_{m}^{M,N}+3Y_{m}^{M,N} ( Z_{m}^{M,N})^2+( Z_{m}^{M,N})^3- Z_{m}^{M,N}}}, AY_{m}^{M,N}\big>
\\
\leq
&
\|Y_{m-1}^{M,N}\|^2
-
 k\|AY_{m}^{M,N}\|^2
-
\tfrac32k \|\nabla( Y_{m}^{M,N})^2\|^2
+
2k \|\nabla Y_{m}^{M,N}\|^2
\\
&
+
C
k(\|Y_{m}^{M,N}\|_{L^4}^4 \|Z_{m}^{M,N}\|_V^2
+
\|Y_{m}^{M,N}\|^2 \|{\color{black}{Z_{m}^{M,N}}}\|_V^4
+
\|{\color{black}{Z_{m}^{M,N}}}\|_{L^6}^6
+
\|{\color{black}{Z_{m}^{M,N}}}\|^2).
\end{split}
\end{equation}
By summation on $m$ and using \eqref{eq:bound-L4-summation}, \eqref{eq:embedding-equatlity-I} and the fact {\color{black}{$Y_0^{M,N}=P_NPX_0$}}, we deduce that for $\kappa=\min\{\gamma,\frac d2+\frac14\}$,
\begin{align}
\begin{split}
\|Y_{m}^{M,N}\|^2
&
+
\sum_{j=1}^{m}
k\|AY_j^{M,N}\|^2
+
\frac32 \sum_{j=1}^{m}
k\|\nabla(Y_j^{M,N})^2\|^2
\\
\leq
&
\|X_0\|^2
+
2k\sum_{j=1}^{m}|Y_j^{M,N}|_1^2
+
C
\Big( 1 + \sup_{1 \leq j \leq M }\|{\color{black}{Z_j^{M,N}}}\|_V^6 \Big)
\big(1+\sum_{j=1}^{m}k\|Y_j^{M,N}\|_{L^4}^4\big)
\\
\leq
&
C\Big( 1+\sup_{ 1 \leq j \leq M }\|{\color{black}{Z_j^{M,N}}}\|_\kappa^6 \Big)
\Big( 1+\sum_{j=1}^{m}k\|Y_j^{M,N}\|_{L^4}^4 \Big).
\end{split}
\end{align}
Therefore, by using \eqref{eq:bound-full-discrete-stochastic-convulution} and \eqref{eq:bound-L4-summation-expectation}, it enables us to obtain
\begin{align}\label{eq:Y-H2-norm-integrand}
\Big(
\mathbb{E}
\big[
\sup_{ 1\leq m\leq M }
\|Y_m^{M,N}\|^{2p}
\big]
\Big)^{\frac1p}
+
\Big\|k\sum_{j=1}^M\|AY_j^{M,N}\|^2\Big\|_{L^p(\Omega;\mathbb{R})}
+
\Big\|k\sum_{j=1}^M\|\nabla( Y_j^{M,N})^2\|^2\Big\|_{L^p(\Omega;\mathbb{R})}
<\infty.
\end{align}

Now we focus on the boundedness of $Y_m^{M,N}$ in the norm $L^6$. Similarly to the continuous case,  we also consider two cases: either $d=1$ or $d=2,3$. {\color{black}{For $d=1$, {\color{black}{we follow the similar skill used in the proof of \cite[Lemma 4.4]{cui2021strong}  to derive}}
\begin{align}\label{eq:Yhn-L6-norm}
\begin{split}
\mathbb{E}
\Big[ \sup_{ 1\leq m\leq M } \| Y_m^{M,N} \|^p_{L^6} \Big]
\leq
 C(p,T,\gamma,X_0).
\end{split}
\end{align}}}
For $d=2,3$, we  introduce the difference operator  $\overline{\partial} Y_m^{M,N}:=\frac{Y_m^{M,N}-Y^{M,N}_{m-1}}k$ and multiply \eqref{eq:full-modified-problem} by $A^{-1}\overline{\partial} Y_m^{M,N}$ to obtain
\begin{align}
\begin{split}
\big|\overline{\partial} Y_m^{M,N}\big|_{-1}^2
&+
\big<A Y_m^{M,N}, \overline{\partial} Y_m^{M,N}\big>
+
\big<F(Y_m^{M,N}), \overline{\partial} Y_m^{M,N}\big>
\\&=
-
\big<{\color{black}{3 (Y_m^{M,N})^2 Z_m^{M,N} +3 Y_m^{M,N} (Z_m^{M,N})^2+ (Z_m^{M,N})^3-Z_m^{M,N}}},\overline{\partial} Y_m^{M,N}\big>.
\end{split}
\end{align}
{\color{black}{Additionally, by using \eqref{eq:definition-lyapunov-function},
one can find that
\begin{align}
J(Y_m^{M,N})
-J(Y_{m-1}^{M,N})
=
\big<A Y_m^{M,N}, \overline{\partial} Y_m^{M,N}\big>
+
\big<F(Y_m^{M,N}), \overline{\partial} Y_m^{M,N}\big>
+
\big<Y_m^{M,N}-Y_{m-1}^{M,N},\overline{\partial} Y_m^{M,N}\big>.
\end{align}}}
{\color{black}{Then, the above estimates together with the fact $\overline{\partial} Y_m^{M,N}\in \dot{H}$ yield}}
\begin{align}
\begin{split}
&k\big|\overline{\partial} Y_m^{M,N}\big|_{-1}^2
+
J(Y_m^{M,N})
-J(Y_{m-1}^{M,N})
\\
&
=-
k\big<{\color{black}{3 (Y_m^{M,N})^2 Z_m^{M,N} +3 Y_m^{M,N}(Z_m^{M,N})^2+ (Z_m^{M,N})^3-Z_m^{M,N}}},\overline{\partial} Y_m^{M,N}\big>
+
k\big<Y_m^{M,N}-Y_{m-1}^{M,N},\overline{\partial} Y_m^{M,N}\big>
\\
&
\leq
k|{\color{black}{P\big(3(Y_m^{M,N})^2 Z_m^{M,N} +3 Y_m^{M,N} (Z_m^{M,N})^2+ (Z_m^{M,N})^3-Z_m^{M,N}}}\big)|_1^2
+
k|Y_m^{M,N}-Y_{m-1}^{M,N}|_1^2
+
\frac k2|\overline{\partial} Y_m^{M,N}|_{-1}^2
\\
&
\leq
C(\gamma)(1+k\|(Y_m^{M,N})^2\|_1^2+
k\|AY_m^{M,N}\|^2)( 1+\|Z_m^{M,N}\|^3_\gamma)^2
+
2k(|Y_m^{M,N}|_1^2+|Y_{m-1}^{M,N}|_1^2)
+
\frac k2|\overline{\partial} Y_m^{M,N}|_{-1}^2
,
\end{split}
\end{align}
where in the last inequality we applied the similar arguments used in \eqref{eq:v-convulution-in-h1-norm}.
Then summation on $m$ and applying \eqref{eq:bound-L4-summation-expectation} and \eqref{eq:Y-H2-norm-integrand} lead to
\begin{align}
\begin{split}
\mathbb{E}
\Big[
\sup_{ 1\leq m\leq M }\big(J(Y_m^{M,N})\big)^p
\Big]
&\leq
C(\gamma,T)
\Big(1+\mathbb{E}\Big[\Big(k\sum_{j=1}^M\|(Y_j^{M,N})^2\|_1^2\Big)^{2p}
\\
&\quad
+
\Big(k\sum_{j=1}^M\|AY_j^{M,N}\|^2\Big)^{2p}\Big]\Big)^{\tfrac12}
\Big( 1+\mathbb{E}\big[\big(\sup_{1 \leq j \leq M}\|{\color{black}{Z_j^{M,N}}}\|^{12p}_\gamma\big)\big]\Big)^{\tfrac12}
\\
&\quad+
\mathbb{E}\Big[
k\sum_{j=1}^M\|AY_j^{M,N}\|^2\Big]^p
+
\mathbb{E}\left[\left(J(P_NPX_0)\right)^p\right]
+
C\mathbb{E}\left[|P_NPX_0|_1^{2p}\right]
\\&
\leq
C(\gamma,p,T)
.
\end{split}
\end{align}
The above estimate together with \eqref{eq:Y-H2-norm-integrand} and the fact $|Y_m^{M,N}|_1^2{\color{black}{-2\|Y_m^{M,N}\|^2}}\leq 2J(Y_m^{M,N})$ yields
\begin{align}
\mathbb{E}
\Big[
\sup_{1 \leq m \leq M}|Y_m^{M,N}|_1^{2p}
\Big]
<
\infty,
\end{align}
which in combination  with \eqref{eq:Yhn-L6-norm} and {\color{black}{the fact $Y_m^{M,N}\in \dot{H}$}} arrives at, for $d=1,2,3$,
\begin{align}
\mathbb{E}\Big[ \sup_{1\leq m\leq M } \|Y_m^{M,N}\|_{L^6}^p \Big]
<
\infty.
\end{align}
Therefore,
 by \eqref{lem:bounded-time-full-deterministic-problem} and \eqref{eq:bound-full-discrete-stochastic-convulution}, we have, for $\kappa=\min\{\gamma,\frac d2+\frac14\}$,
\begin{align}
\begin{split}
\|Y_m^{M,N}\|_{L^p(\Omega;H^\kappa)}
\leq
&
\|E_{k,N}^mPX_0\|_{L^p(\Omega;H^\kappa)}
+
\big\|
\sum_{i=1}^{m}kE_{k,N}^{m-i+1}APF(Y_i^{M,N}+{\color{black}{Z_i^{M,N}}})
\big\|_{L^p(\Omega;H^\kappa)}
\\
\leq
&
\|X_0\|_{L^p(\Omega;H^\kappa)}
+
\sum_{i=1}^{m}
kt_{m-i+1}^{-\frac{2+\kappa}4}
\|F(Y_i^{M,N}+{\color{black}{Z_i^{M,N}}})
\|_{L^p(\Omega;H)}
\\
\leq
&
\|X_0\|_{L^p(\Omega;H^\gamma)}
+
\sum_{i=1}^{m}kt_{m-i+1}^{-\frac{2+\kappa}4}\sup_{1 \leq i \leq M}\big(\|Y_i^{M,N}\|^3_{L^{3p}(\Omega;L^6)}+
\|{\color{black}{Z_i^{M,N}}}\|^3_{L^{3p}(\Omega;L^6)}+1\big)
\\
\leq
&C(p,T,\gamma,X_0).
\end{split}
\end{align}
Hence this finishes the proof of this theorem. $\square$
\section{Strong convergence rates of the full discretization}
In this part, we follow the approach in \cite{ruisheng2018error} to derive the error estimates of the fully discrete  approximation of the stochastic problem \eqref{eq:abstract-CHC-form}.
The convergence analysis heavily relies on the moment bound obtained in section 4 and the corresponding deterministic error estimates.

The next theorem states the main result of this paper, concerning strong convergence rates of the  full discretizaiton scheme.
\begin{theorem}\label{them:error-estimates-full-stochastic}
Let $X(t)$ be the mild solution of \eqref{eq:abstract-CHC-form} and  $X_m^{M,N}$ be the solution of \eqref{eq:full-discretization}. Suppose
Assumptions \ref{assum:linear-operator-A}-\ref{assum:intial-value-data} are valid, then
 there exist two positive constants $C$ and $k_0$ such that for all $k\leq k_0$ and $\forall p\geq 1$,
\begin{align}\label{them:error-estimates-full-problem}
\sup_{1 \leq m \leq M }
\|X(t_m)-X_m^{M,N}\|_{L^p(\Omega;H)}
\leq
C(\lambda_N^{-\frac\gamma2} +k^{\frac\gamma 4} ).
\end{align}
\end{theorem}
To show this theorem,  we introduce an auxiliary problem,
\begin{align}\label{eq:auxiliary-problem}
\widetilde{X}_m^{M,N}-\widetilde{X}_{m-1}^{M,N}
+
k A^2 \widetilde{X}_m^{M,N}
+
kP_NAF(X(t_m))
=
P_N \Delta W_m,\; X_0^{M,N}={\color{black}{P_N}}X_0,
\end{align}
whose solution can be recasted as
\begin{align}\label{eq:solution-form-auxiliary-problem}
{\color{black}{\widetilde{X}_m^{M,N}={\color{black}{E_{k,N}^mP_N}}X_0-k\sum_{j=1}^mE_{k,N}^{m-j+1}APF(X(t_j))
+Z_m^{M,N}}}.
\end{align}
Owning to  \eqref{lem:bounded-time-full-deterministic-problem}, \eqref{lem:eq-bound-sum-full-operator-II} and \eqref{eq:bound-full-discrete-stochastic-convulution}, one can derive that, for $\kappa=\min\{\gamma,\frac d2+\frac14\}$ and any $m\in\{1,2,\cdots,M\}$,
\begin{align}\label{eq:boundness-solution-form-full-auxiliary-problem}
\begin{split}
\|\widetilde{X}_m^{M,N}\|_{L^p(\Omega;H^\kappa)}
\leq&
\|{\color{black}{E_{k,N}^mP_NX_0}}\|_{L^p(\Omega;H^\kappa)}
+
k\sum_{j=1}^m\|E_{k,N}^{m-j+1}APF(X(t_j))\|_{L^p(\Omega;H^\kappa)}
+
\|{\color{black}{Z_m^{M,N}}}\|_{L^p(\Omega;H^\kappa)}
\\
\leq
&
C\|X_0\|_{L^p(\Omega;H^\kappa)}
+
Ck\sum_{j=1}^mt_{m-j+1}^{-\frac{2+\kappa}4}\|F(X(t_j))\|_{L^p(\Omega;H)}
+
C\|A^{\frac{\gamma-2}2}Q^{\frac12}
\|_{\widetilde{\mathcal{L}}_2(H)}
\\
\leq
&
C\big(1+\|X_0\|_{L^p(\Omega;H^\gamma)}
+
k\sum_{j=1}^mt_{m-j+1}^{-\frac{2+\kappa}4}
\sup_{s\in[0,T]}\|F(X(s))\|_{L^p(\Omega;H)}
\big)
<\infty.
\end{split}
\end{align}
Hence, we can decompose the considered error $\|X(t_m)-X_m^{M,N}\|_{L^p(\Omega;H)}$ into two parts:
\begin{align}\label{eq:error-decompose}
\|X(t_m)-X_m^{M,N}\|_{L^p(\Omega;H)}
\leq
\|X(t_m)-\widetilde{X}_m^{M,N}\|_{L^p(\Omega;H)}
+
\|\widetilde{X}_m^{M,N}-X_m^{M,N}\|_{L^p(\Omega;H)}.
\end{align}
These two error terms are separately handled in the forthcoming two subsections.
Plugging \eqref{them:error-estimates-auxiliary-problem} and \eqref{them:error-estimates-auxiliary-problem}
into \eqref{eq:error-decompose} gives the desired assertion \eqref{them:error-estimates-full-problem}.
%
\subsection{Error estimates for $\|X(t_m)-\widetilde{X}_m^{M,N}\|_{L^p(\Omega;H)}$}
In what follows, we use the corresponding deterministic error estimates to bound $\|X(t_m)-\widetilde{X}_m^{M,N}\|_{L^p(\Omega;\dot{H})}$
in the semigroup framework, which is usually applied in the error analysis of numerical approximation of  SPDEs
with globally Lipschitz coefficients.
%
\begin{lemma}\label{lem:bound-first-error-term}
Let $X(t)$ be the mild solution of \eqref{eq:abstract-CHC-form} and  $\widetilde{X}_m^{M,N}$ be the solution of \eqref{eq:auxiliary-problem}. Suppose
Assumptions \ref{assum:linear-operator-A}-\ref{assum:intial-value-data} are valid,
then there exist two positive constants $C$ and $k_0$ such that for all $k\leq k_0$ and $\forall p\geq 1$,
\begin{align}\label{them:error-estimates-auxiliary-problem}
\|X(t_m)-\widetilde{X}_m^{M,N}\|_{L^p(\Omega;H)}
\leq
C(\lambda_N^{-\frac\gamma2} +k^{\frac\gamma 4} ).
\end{align}
\end{lemma}
Its proof is {\color{black}{given below.}}
First, we define the fully discrete approximation operator
\begin{align}
\Psi_k^{M,N}(t):
=
E(t)-E_{k,N}^m,\;\;\;t\in[t_{m-1},t_m),\;m\in\{1,2,\cdots,M\}.
\end{align}
In the following lemma, we give some results on   the error operator $\Psi_k^{M,N}(t)$, which  are crucial in the error estimates of the  fully discrete approximation and can be shown by a slight
 modification of the proof of  \cite[Lemma 6.2]{ruisheng2018error}.
\begin{lemma}\label{lem:error-estimes-time-full-deterministic-problem}
Under Assumption  \ref{assum:linear-operator-A}, the following estimates for $\Psi^{M,N}_k(t)$ hold for $t\in[0,T]$.
\\
(i) Let $\beta\in[0,4]$. There exists a constant $C$ such that, for $t>0$,
\begin{align}\label{lem:error-stimates-nonsmooth}
\|\Psi_k^{M,N}(t) v\|
&\leq
C(\lambda_N^{-\frac\beta2}+k^{\frac \beta 4})  |v|_\beta,
\quad
\forall x\in {\color{black}{H^\beta}}.
\end{align}
(ii) Let $ \alpha\in[0,2]$. There exists a constant $C$ such that, for $t>0$,
\begin{align}
\|\Psi_k^{M,N}(t) v\|
&\leq
C(\lambda_N^{-\frac{4-\alpha}2}+k^{\frac {4-\alpha} 4}) t^{-1} |v|_{-\alpha},
\quad
\forall v\in {\color{black}{H^{-\alpha}}}.
\label{lem:eq-error-error-estimes-time-full-deterministic-problem}
\end{align}
(iii)
Let $\nu\in[0,4]$. There exists a constant $C$ such that, for $t>0$,
\begin{align}\label{lem:error-F-full-integrand}
\Big(\int_0^t\|\Psi_k^{M,N}(s)v\|^2\, \mathrm{d} s\Big)^{\frac12}
\leq
C(\lambda_N^{-\frac\nu2}+  k^{\frac\nu4})|v|_{\nu-2},
\quad
v\in {\color{black}{H^{\nu-2}}}.
\end{align}
(v)
Let $\varrho\in[0,2]$. There exists a  constant $C$ such that, for $t>0$,
\begin{align}\label{lem:error-deterministc-potial-full-integrand-II}
\Big\|\int_0^t \Psi_k^{M,N}(s)v \,\mathrm{d} s\Big\|
\leq
C(\lambda_N^{-\frac{4-\varrho}2}+k^{\frac{4-\varrho}4})|v|_{-\varrho},
\quad
v\in {\color{black}{H^{-\varrho}}}.
\end{align}
\end{lemma}
Subsequently, we are well-prepared to show Lemma \ref{lem:bound-first-error-term}.

{\it Proof of Lemma \ref{lem:bound-first-error-term}.}
Subtracting \eqref{eq:solution-form-auxiliary-problem} from
  \eqref{them:eq-mild-solution-stochastic-equation}, the error $X(t_m)-\widetilde{X}_m^{M,N}$ reads
\begin{align}\label{eq:full-auxiliary-error}
\begin{split}
\|X(t_m)-&\widetilde{X}_m^{M,N}\|_{L^p(\Omega;H)}
=
\|(E(t_m)-E_{k,N}^m)X_0\|_{L^p(\Omega;H)}
\\
&+
\Big\|\int_0^{t_m}E(t_m-s)APF(X(s))\,\dd s
-
k\sum_{j=1}^mE_{k,N}^{m-j+1}APF(X(t_j))\Big\|_{L^p(\Omega;H)}
\\
&+
\Big\|\sum_{j=1}^m\int_{t_{j-1}}^{t_j}\big(E(t_m-s)
-
E_{k,N}^{m-j+1}\big) \,\dd W(s)\Big\|_{L^p(\Omega;H)}
\\
=:&
I+J+K.
\end{split}
\end{align}
Subsequently,  $I, J, K$ will be treated separately.
For the first term $I$, we employ  \eqref{lem:error-stimates-nonsmooth} with $\beta=\gamma$ to obtain
\begin{align}
I
\leq
C(\lambda_N^{-\frac\gamma2}+k^{\frac\gamma 4})\|X_0\|_{L^p(\Omega;H^\gamma)}.
\end{align}
To handle $J$,  we decompose it into three terms as follows
\begin{align}
\begin{split}
J
\leq
&
\Big\|
\sum_{j=1}^m\int_{t_{j-1}}^{t_j}E(t_m-s)AP\big(F(X(s))-F(X(t_j))\big)
\,\dd s
\Big\|_{L^p(\Omega;H)}
\\
&
+
\Big\|\int_0^{t_m}\Psi_k^{M,N}(t_m-s)APF(X(t_m))\,\dd s\Big\|_{L^p(\Omega;H)}
\\
&+
\sum_{j=1}^
m\int_{t_{j-1}}^{t_j}\|\Psi_k^{M,N}(t_m-s)AP
\big(F(X(t_j))-F(X(t_m))\big)\|_{L^p(\Omega;H)}\,\dd s
\\
=:
&
J_1+J_2+J_3.
\end{split}
\end{align}
Next we treat the above three terms separately.
To deal with $J_1$, we first  note that, for $s\in[t_{j-1},t_j)$,
\begin{align}
X(t_j)
=
E(t_j-s)X(s)
-
\int_{s}^{t_j}E(t_j-\sigma)AF(X(\sigma))\,\dd \sigma
+
\int_{s}^{t_j}E(t_j-\sigma)\,\dd W(\sigma),
\end{align}
and then apply  Taylor's formula to split $J_1$ into four terms:
\begin{align}
\begin{split}
J_1
\leq
&
\Big\|\sum_{j=1}^m\int_{t_{j-1}}^{t_j}E(t_m-s)APF'(X(s))(E(t_j-s)-I)X(s)\,\dd s\Big\|_{L^p(\Omega;H)}
\\
&+
\Big\|\sum_{j=1}^m\int_{t_{j-1}}^{t_j}E(t_m-s)APF'(X(s))\int_s^{t_j}
E(t_j-\sigma)AF(X(\sigma))\,\dd \sigma\,\dd s\Big\|_{L^p(\Omega;H)}
\\
&+
\Big\|\sum_{j=1}^m\int_{t_{j-1}}^{t_j}E(t_m-s)APF'(X(s))\int_s^{t_j}
E(t_j-\sigma)\,\dd W(\sigma)\,\dd s \Big\|_{L^p(\Omega;H)}
\\
&+
\Big\|\sum_{j=1}^m\int_{t_{j-1}}^{t_j}E(t_m-s)APR_F(X(s),X(t_j))\,\dd s\Big\|_{L^p(\Omega;H)}
\\
=:
&
J_{11}+J_{12}+J_{13}+J_{14},
\end{split}
\end{align}
where $R_F(X(s),X(t_j))$ is a remainder term, given by
\[
R_F(X(s),X(t_j))=\int_0^1
F''\big(X(s)+\lambda(X(t_j)-X(s))\big)(X(t_j)-X(s),X(t_j)-X(s))(1-\lambda)\,\dd \lambda.\]
Owing to \eqref{them:spatial-regularity-mild-stoch}, \eqref{eq:embedding-equatlity-III} and \eqref{I-spatio-temporal-S(t)} with $\mu=\frac{2+\delta_0}2$,
  we deduce, for any fixed $\delta_0\in(\frac32,2)$ and $\gamma \in(\frac d2,4]$,
\begin{align}\label{eq:estimate-L2111}
\begin{split}
J_{11}
&\leq
C\sum_{j=1}^m\int_{t_{j-1}}^{t_j}(t_m-s)^{-\frac{2+\delta_0}4}
      \|A^{-\frac{\delta_0}2}PF'(X(s))(E(t_j-s)-I)X(s)\|_{L^{p}(\Omega;H)}\,\dd s
\\
&\leq
C\sum_{j=1}^m\int_{t_{j-1}}^{t_j}(t_m-s)^{-\frac{2+\delta_0}4}\|F'(X(s))
    (E(t_j-s)-I)X(s)\|_{L^{p}(\Omega;L^1)}\,\dd s
\\
&\leq
C\sum_{j=1}^m
\int_{t_{j-1}}^{t_j}
(t_m-s)^{-\frac{2+\delta_0}4}
\Big(
1+\|X(s)\|^2_{L^{4p}(\Omega;L^4)}
\Big)
\Big\|
(E(t_j-s)-I)X(s)
\Big\|_{L^{2p}(\Omega;H)}\,\dd s
\\
&\leq
C
k^{\frac\gamma4}
\int_0^{t_m}
(t_m-s)^{-\frac{2+\delta_0}4}
\,\dd s
\Big(
1+\sup_{s\in[0,T]}
\|X(s)\|^3_{L^{4p}(\Omega;H^\gamma)}
\Big)
\\
&\leq
 Ck^{\frac\gamma4}.
\end{split}
\end{align}
Using the similar approach as used in the proof of \eqref{eq:estimate-L2111}, with \eqref{lem:eq-F-mild-bound} used instead we obtain,
\begin{align}\label{eq:estimate-L2112}
\begin{split}
J_{12}
&\leq C
\sum_{j=1}^m\int_{t_{j-1}}^{t_j}\int_{s}^{t_j}(t_m-s)^{-\frac{2+\delta_0}4}
\big\|F'(X(s))E(t_{j}-\sigma)AF(X(\sigma))\big\|_{L^{p}(\Omega;L^1)}\,\dd \sigma\,\dd s
\\
&\leq
C\sum_{j=1}^{m}\int_{t_{j-1}}^{t_j}\int_s^{t_j}(t_m-s)^{-\frac{2+\delta_0}4}
(t_{j}-\sigma)^{-\frac{2-{\ell(\gamma)}}4}
(1+\|X(s)\|^2_{L^{4p}(\Omega;L^4)})
    \|F(X(\sigma))\|_{L^{2p}(\Omega;H^{\ell(\gamma)})}
    \,\dd \sigma\,\dd s
\\
&\leq
Ck^{\frac{2+\ell(\gamma)}4}
\int_0^{t_m}(t_m-s)^{-\frac{2+\delta_0}4}
\dd s
\Big(
1+\sup_{s\in[0,T]}\|X(s)\|^2_{L^{4p}(\Omega;H^{\gamma})}
\Big)
\sup_{s\in[0,T]}\|F(X(s))\|_{L^{2p}(\Omega;H^{\ell(\gamma)})}
\\
&\leq
 Ck^{\frac\gamma4},
\end{split}
\end{align}
where
$\ell(\gamma)$ is a piecewise function defined by $\ell(\gamma)
=
0
\text{ for } \gamma\in(\frac d2,2]
\text{ and }
\ell(\gamma) = 2
\text{ for } \gamma\in (2,4].
$
To bound $J_{13}$, we introduce an indicator function defined by
$\chi_{[s,t_j)}(t)
=
1
\text{ for } t\in [s,t_j)
\text{ and }
\chi_{[s,t_j)}(t) = 0
\text{ for } t\not \in [s,t_j)
$ and then employ the stochastic Fubini theorem (see \cite[Theorem 4.18]{da2014stochastic}) and the Burkholder-Davis-Gundy-type inequality to obtain
\begin{align}
\begin{split}
J_{13}
=&
\Big\|\sum_{j=1}^{m}\int_{t_{j-1}}^{t_j}\int_{t_{j-1}}^{t_j}\chi_{[s,t_j)}(\sigma)
   E(t_m-s)APF'(X(s))E(t_j-\sigma)\,\dd W(\sigma)\,\dd s \Big\|_{L^{p}(\Omega;H)}
\\
=
&
\Big\|\sum_{j=1}^{m}\int_{t_{j-1}}^{t_j}\int_{t_{j-1}}^{t_j}\chi_{[s,t_j)}(\sigma)
  E(t_m-s)APF'(X(s))E(t_j-\sigma)\,\dd s \dd W(\sigma) \Big\|_{L^{p}(\Omega;H)}
\\
\leq
&
\Big(\sum_{j=1}^{m}\int_{t_{j-1}}^{t_j}\Big\|\int_{t_{j-1}}^{t_j}
\chi_{[s,t_j)}(\sigma)E(t_m-s)APF'(X(s))E(t_j-\sigma)Q^{\frac12}\,\dd s\Big\|_{L^p(\Omega;\mathcal{L}_2(H))}^2\,\dd \sigma\Big)^{\frac12}.
\end{split}
\end{align}
Further, employing \eqref{them:spatial-regularity-mild-stoch}, \eqref{eq:FX1-FX2-d=23},  \eqref{eq:FX1-FX2-d=1} with {\color{black}{$\iota=\kappa=\min\{\gamma,\frac d2+\frac14\}$}},  \eqref{I-spatio-temporal-S(t)} with  $\mu=\frac{2+\eta}2$ and \eqref{III-spatio-temporal-S(t)} with $\varrho=\frac{\max\{3-\gamma,0\}}2$ and using the H\"{o}lder inequality, one can find that
\begin{align}\label{eq:j13-estimatation-full}
\begin{split}
J_{13}
\leq
&
C
k^{\frac12}
\Big(
\sum_{j=1}^m
\int_{t_{j-1}}^{t_j}
\int_{t_{j-1}}^{t_j}
\sum_{l=1}^\infty
\|
E(t_m-s)APF'(X(s))E(t_j-\sigma)Q^{\frac12}\eta_l
\|_{L^p(\Omega;H)}^2
\,\dd s\,\dd \sigma
\Big)^{\frac12}
\\
\leq
&
Ck^{\frac12}
 \Big(
 \sum_{j=1}^m
 \int_{t_{j-1}}^{t_j}
\int_{t_{j-1}}^{t_j}
(t_m-s)^{-\frac{2-\eta}2}
\sum_{l=1}^\infty{\color{black}{\|
F'(X(s))E(t_j-\sigma)Q^{\frac12}\eta_l
\|_{L^p(\Omega;H^\eta)}^2}}
\,\dd s\dd \sigma
\Big)^{\frac12}
\\
\leq
&
Ck^{\frac12}
\Big(
\sum_{j=1}^m
\int_{t_{j-1}}^{t_j}
(t_m-s)^{-\frac{2-\eta}2}
\Big(
1+\| X(s)\|^4_{L^{2p}(\Omega;H^\kappa)}
\Big)\,\dd s
\int_{t_{j-1}}^{t_j}
\sum_{l=1}^\infty
{\color{black}{\|
E(t_j-\sigma)Q^{\frac12}\eta_l
\|_1^2}}
\dd \sigma
\Big)^{\frac12}
\\
\leq
&
{\color{black}{
Ck^{\frac12}
\Big(
\int_0^{t_m}
(t_m-s)^{-\frac{2-\eta}2}\,\dd s\Big)^{\frac12}
\Big(
1+
\sup_{s\in[0,T]}
\| X(s)\|^4_{L^{2p}(\Omega;H^\gamma)}
\Big)^{\frac12}}}
\\
&
\qquad
\qquad
{\color{black}{
\sup_{j\in\{1,2,\cdots,M\}}\Big(
\sum_{l=1}^\infty
\int_{t_{j-1}}^{t_j}\big(\|A^{\frac12}E(t_j-\sigma)P Q^{\frac12}\eta_l\|^2
+
\big< Q^{\frac12}\eta_l,e_0\big>^2\big)\,\dd \sigma\Big)^{\frac12}}}
\\
\leq
&
{\color{black}{
C
k^{\frac{4-\max\{3-\gamma,0\}}4}
\big(\sum_{l=1}^\infty\|A^{\frac{\gamma-2}2}P Q^{\frac12}\eta_l\|^2+\big< Q^{\frac12}\eta_l,e_0\big>^2\big)^{\frac12}}}
\\
\leq
&
{\color{black}{
Ck^{\frac\gamma4}(\|A^{\frac{\gamma-2}2}PQ^{\frac12}\|_{\mathcal{L}_2(H)}+\|(I-P)Q^{\frac12}\|_{\mathcal{L}_2(H)}),
}}
\end{split}
\end{align}
where  $\eta := \min\{\gamma,\frac23\}$,   and $\{\eta_l\}_{l\in\N^+}$ is the eigenfunction of $Q$.
Now we are in the position to handle the term $J_{14}$.
It follows from \eqref{them:temporal-regularity-mild-stoch} and \eqref{them:spatial-regularity-mild-stoch} that
\begin{align}
\begin{split}
J_{14}
&\leq
C
\sum_{j=1}^m
\int_{t_{j-1}}^{t_j}\!\!
(t_m-s)^{-\frac{2+\delta_0}4}
\\
&
\qquad
\Big\|
\int_0^1 \!\!
F''\big(X(s)+\lambda(X(t_j)-X(s))\big)(X(t_j)-X(s),X(t_j)-X(s))
\,\dd \lambda
\Big\|_{L^p(\Omega;L^1)}
\,\dd s
\\
&
\leq
C
\sum_{j=1}^m
\int_{t_{j-1}}^{t_j}
(t_m-s)^{-\frac{2+\delta_0}4}
\|X(t_j)-X(s)\|^2_{L^{4p}(\Omega;H)}
\,
\big(
1+\sup_{s\in[0,T]}\|X(s)\|_{L^{2p}(\Omega;V)}
\big)
\,\dd s
\\
&\leq
Ck^{\min\{1,\frac\gamma2\}}
\Big(
1+\sup_{s\in[0,T]}\|X(s)\|_{L^{2p}(\Omega;H^\gamma)}
\Big)
         \int_0^{t_m}(t_m-s)^{-\frac{2+\delta_0}4}\,\dd s
\\
&
\leq Ck^{\frac\gamma4},
\end{split}
\end{align}
where in the first inequality we used similar arguments as in \eqref{eq:estimate-L2111}. Thus, putting the above four estimates together results in, for $\gamma\in(\frac d2,4]$,
\begin{align}\label{eq:bound-J1}
J_1
\leq
Ck^{\frac\gamma4}.
\end{align}
Concerning  the term $J_2$, by  applying  \eqref{lem:eq-F-mild-bound} and \eqref{lem:error-deterministc-potial-full-integrand-II} with $\varrho=2-\ell(\gamma)$, one can observe that
\begin{align}\label{eq:bound-J2}
\begin{split}
J_2
&\leq
C(\lambda_N^{-\frac{2+\ell(\gamma)}2}+k^{\frac{2+\ell(\gamma)}4})
\|F(X(t_m))\|_{L^p(\Omega;H^{\ell(\gamma)})}
\\
&\leq
C(\lambda_N^{-\frac{2+\ell(\gamma)}2} +k^{\frac{2+\ell(\gamma)}4}) \sup_{s\in[0,T]} \|F(X(s))\|_{L^p(\Omega;H^{\ell(\gamma)})}
\leq
C(\lambda_N^{-\frac\gamma2}+k^{\frac\gamma4}),
\end{split}
\end{align}
where $\ell(\gamma)=0$, for $\gamma\in(\tfrac d2,2]$ and $\ell(\gamma)=2$, for $\gamma\in(2,4]$ by \eqref{eq:l-gamma-defn}.
With regard to $J_3$, we  employ \eqref{lem:eq-error-error-estimes-time-full-deterministic-problem} with $\alpha=4-l(\gamma)$
 and \eqref{lem:temporal-regulairty-F}  to infer
\begin{align}
\begin{split}
J_3
&\leq
C(\lambda_N^{-\frac{l(\gamma)}2} + k^{\frac{l(\gamma)} 4} )\sum_{j=1}^
m\int_{t_{j-1}}^{t_j}(t_m-s)^{-1}\|F(X(t_j)-F(X(t_m))
\|_{L^p(\Omega;H^{l(\gamma)-2})}\,\dd s
\\
&\leq
C(\lambda_N^{-\frac{l(\gamma)}2} + k^{\frac{l(\gamma)} 4} )\sum_{j=1}^
{m-1}\int_{t_{j-1}}^{t_j}(t_m-s)^{-1} t_{m-j}^{\frac18}
\,\dd s
\\
&\leq
C(\lambda_N^{-\frac\gamma2}+ k^{\frac\gamma 4} ){\color{black}{\int_0^{t_{m-1}}(t_m-s)^{-1+\frac18}\,\dd s}}
\\
&
\leq C(\lambda_N^{-\frac\gamma2}+ k^{\frac\gamma 4} ),
\end{split}
\end{align}
{\color{black}{where
$l(\gamma)=2$, for $\gamma\in(\tfrac d2,2]$, $l(\gamma)=3$, for $\gamma\in(2,3]$ and $l(\gamma)=4$, for $\gamma\in(3,4]$ by \eqref{eq:l-gamma-defn}.}}
This combined with \eqref{eq:bound-J1} and \eqref{eq:bound-J2} leads to
\begin{align}
J
\leq
C(\lambda_N^{-\frac\gamma2} + k^{\frac \gamma 4} ).
\end{align}
To bound  $K$,
utilizing the Burkholder-Davis-Gundy-type inequality and  \eqref{lem:error-F-full-integrand}
with $\nu=\gamma$ arrives at
\begin{align}
\begin{split}
K
\leq
&
C\Big(\sum_{j=1}^m\int_{t_{j-1}}^{t_j}\|\Psi_k^{M,N}(t_m-s)Q^{\frac12}
\|^2_{\mathcal{L}_2(H)}\dd s\Big)^{\frac12}
\\
\leq
&
C(\lambda_N^{-\frac\gamma2}+k^{\frac\gamma 4} ) {\color{black}{(\|A^{\frac{\gamma-2}2}PQ^{\frac12}\|_{\mathcal{L}_2(H)}+\|(I-P)Q^{\frac12}\|_{\mathcal{L}_2(H)}).}}
\end{split}
\end{align}
At last, gathering  the  estimates of $I$, $J$ and $K$ together implies
\begin{align}\label{eq:x-x^{M,N}-estiamte}
\|X(t_m)-\widetilde{X}_m^{M,N}\|_{L^p(\Omega;\dot{H})}
\leq
C(\lambda_N^{-\frac{\gamma}2}+k^{\frac\gamma4}).
\end{align}
Hence this ends the proof of this lemma. $\square$

\subsection{Error estimates for $\|\widetilde{X}_m^{M,N}-X_m^{M,N}\|_{L^p(\Omega;H)}$}
This part is devoted to the bound of the second error term in \eqref{eq:error-decompose}.
The bound of the first error term together with  the moment bounds of the fully discrete approximation and the solution of the auxiliary problem \eqref{eq:auxiliary-problem}  plays a key role in our convergence analysis.
\begin{lemma}\label{lem:bound-second-error-term}
Let  $X_m^{M,N}$ and $\widetilde{X}_m^{M,N}$ be the solutions of \eqref{eq:full-discretization} and \eqref{eq:auxiliary-problem}, respectively. Suppose
Assumptions \ref{assum:linear-operator-A}-\ref{assum:intial-value-data} are valid, then
 there exist two positive constants $C$ and $k_0$ such that for all $k\leq k_0$ and $\forall p\geq 1$,
\begin{align}\label{them:error-estimates-auxiliary-problem}
\|X_m^{M,N}-\widetilde{X}_m^{M,N}\|_{L^p(\Omega;H)}
\leq
C(\lambda_N^{-\frac\gamma2} +k^{\frac\gamma 4} ).
\end{align}
\end{lemma}
{\it Proof of Lemma \ref{lem:bound-second-error-term}.}
{\color{black}{Note first that the  error $\widetilde{e}_m^{M,N}:=X_m^{M,N}-\widetilde{X}_m^{M,N}\in \dot{H}$ is}} the solution of  the following error problem
\begin{align}\label{eq:error-equaton-discrete}
\widetilde{e}_m^{M,N}-\widetilde{e}_{m-1}^{M,N}
+
kA^2\widetilde{e}_m^{M,N}=-kP_NAPF(X_m^{M,N})+kP_NAPF(X(t_m)),\quad \widetilde{e}_0^{M,N}=0,
\end{align}
which can be reformulated as
\begin{align}\label{eq:solution-error-equation-full-discrete}
\widetilde{e}_m^{M,N}=k\sum_{j=1}^mE_{k,N}^{m-j+1}P_NAP(F(X(t_j))-F(X_j^{M,N})).
\end{align}
Multiplying both sides of \eqref{eq:error-equaton-discrete} by $A^{-1}\widetilde{e}_m^{M,N}$ and using \eqref{eq:local-condition}  yield
\begin{align}
\begin{split}
\big<
\widetilde{e}_m^{M,N}-&\widetilde{e}_{m-1}^{M,N},
A^{-1}\widetilde{e}_m^{M,N}
\big>
+
k|\widetilde{e}_m^{M,N}|_1^2
\\
&=
k\big<F(\widetilde{X}_m^{M,N})-F(X_m^{M,N}),\widetilde{e}_m^{M,N}\big>
+
k\big<F(X(t_m))-F(\widetilde{X}_m^{M,N}),\widetilde{e}_m^{M,N}\big>
\\
&\leq
\tfrac{3k}2\|\widetilde{e}_m^{M,N}\|^2
+
\tfrac k2\|F(X(t_m))-F(\widetilde{X}_m^{M,N})\|^2
\\
&\leq
\tfrac{3k}2\|\widetilde{e}_m^{M,N}\|^2
+
Ck\|X(t_m)-\widetilde{X}_m^{M,N}\|^2(1+\|X(t_m)\|^4_V+\|\widetilde{X}_m^{M,N}\|^4_V).
\end{split}
\end{align}
Noting that $\|\widetilde{e}_m^{M,N}\|^2\leq
|\widetilde{e}_m^{M,N}|_1|\widetilde{e}_m^{M,N}|_{-1} $ and $\frac12(|\widetilde{e}_m^{M,N}|_{-1}^2-|\widetilde{e}^{M,N}_{m-1}|_{-1}^2)
\leq
\big<\widetilde{e}_m^{M,N}-\widetilde{e}^{M,N}_{m-1},A^{-1}\widetilde{e}_m^{M,N}\big>$,
we further arrive at
\begin{align}
\begin{split}
\tfrac12(|\widetilde{e}_m^{M,N}|_{-1}^2&-|\widetilde{e}_{m-1}^{M,N}|_{-1}^2)
+
k| \widetilde{e}_m^{M,N}|_1^2
\\
&
\leq
\big<\widetilde{e}_m^{M,N}-\widetilde{e}^{M,N}_{m-1},A^{-1}\widetilde{e}_m^{M,N}\big>
+
k|\widetilde{e}_m^{M,N}|_1^2
\\
&
\leq
\tfrac{k}2|\widetilde{e}_m^{M,N}|^2_1
+
\tfrac{9k}8|\widetilde{e}_m^{M,N}|^2_{-1}
+
Ck\|X(t_m)-\widetilde{X}_m^{M,N}\|^2\big(1+\|X(t_m)\|^4_V+
\|\widetilde{X}_m^{M,N}\|^4_V\big).
\end{split}
\end{align}
By iteration and using the {\color{black}{Gronwall}} inequality, we have
\begin{align}
|\widetilde{e}_m^{M,N}|_{-1}^2
+
k\sum_{j=1}^m|\widetilde{e}^{M,N}_j|_1^2
\leq
Ck\sum_{j=1}^m \|\widetilde{X}_j^{M,N}-X(t_j)\|^2\big(1+\|\widetilde{X}_j^{M,N}\|_V^4
+\|X(t_j)\|_V^4\big).
\end{align}
Further, taking  \eqref{them:spatial-regularity-mild-stoch}, \eqref{eq:boundness-solution-form-full-auxiliary-problem},
\eqref{eq:embedding-equatlity-I} and \eqref{eq:x-x^{M,N}-estiamte} into account,
it follows that, for $\kappa=\min\{\gamma,\frac d2+\frac14\}$
\begin{align}\label{eq:sum-x-x^{M,N}-in-H1}
\begin{split}
\Big\|k\sum_{j=1}^{m}|\widetilde{e}_j^{M,N}|_1^2\Big\|_{L^p(\Omega;\mathbb{R})}
&\leq
C
k\sum_{j=1}^{m} \Big\|\|\widetilde{X}_j^{M,N}-X(t_j)\|^2
\big(1+\|\widetilde{X}_j^{M,N}\|_V^4+\|X(t_j)\|_V^4\big)
\Big\|_{L^p(\Omega;\mathbb{R})}
\\
&\leq
Ck\sum_{j=1}^{m} \|\widetilde{X}_j^{M,N}-X(t_j)\|_{L^{4p}(\Omega; H)}^2
\big(1+\|\widetilde{X}_j^{M,N}\|_{L^{8p}(\Omega;H^\kappa)}^4
+\|X(t_j)\|_{L^{8p}(\Omega;H^\kappa)}^4\big)
\\
&\leq
C(\lambda_N^{-\frac\gamma2} +k^{\frac\gamma 4})^2.
\end{split}
\end{align}
In analogy with \eqref{eq:decompose-en}, we  split $\big\|\widetilde{e}_m^{M,N}\big\|_{L^p(\Omega;H)}$
into two terms:
\begin{align}\label{eq:widetilde{e}-decompose-full}
\begin{split}
\|\widetilde{e}^{M,N}_m\|_{L^p(\Omega;H)}
\leq&
k\sum_{j=1}^m
\Big\|
E_{k,N}^{m-j+1}P_N A P
\big(
F(X(t_j))-F(\widetilde{X}^{M,N}_j)
\big)
\Big\|_{L^p(\Omega;H)}
\\
&+
k\Big\|
\sum_{j=1}^m
E_{k,N}^{m-j+1}P_NAP
\big(
F(\widetilde{X}_j^{M,N})-F(X_j^{M,N})
\big)
\Big\|_{L^p(\Omega;H)}
\\
=:
&
L_1+L_2.
\end{split}
\end{align}
Following similar arguments used in the proof of \eqref{eq:sum-x-x^{M,N}-in-H1}, we obtain, for $\kappa=\min\{\gamma,\frac d2+\frac14\}$
\begin{align}\label{eq:A-estimate-full}
\begin{split}
L_1
&\leq
C
k
\sum_{j=1}^m
t_{m-j+1}^{-\frac12}
\|F(X(t_j))-F(\widetilde{X}_j^{M,N})\|_{L^p(\Omega;H)}
\\
&
\leq C
k
\sum_{j=1}^m
t_{m-j+1}^{-\frac12}
\|
X(t_j)-\widetilde{X}_j^{M,N}
\|_{L^{2p}(\Omega;H)}
\Big(
1+\|X(t_j)\|^2_{L^{4p}(\Omega;V)}
+
\|\widetilde{X}_j^{M,N}\|^2_{L^{4p}(\Omega;V)}
\Big)
\\
&
\leq
C(\lambda_N^{-\frac\gamma2} +k^{\frac\gamma 4})
k
\sum_{j=1}^m
t_{m-j+1}^{-\frac12}
\Big(
1+
\sup_{s\in[0,T]}
\|X(s)\|^2_{L^{4p}(\Omega;H^\kappa)}
+
\sup_{1\leq j\leq M}
\|\widetilde{X}_j^{M,N}\|^2_{L^{4p}(\Omega;H^\kappa)}
\Big)
\\
&
\leq
C(\lambda_N^{-\frac\gamma2} + k^{\frac\gamma 4}).
\end{split}
\end{align}
{\color{black}{To bound  the term $L_2$, we  apply  \eqref{eq:FX1-FX2-d=1} with $\iota=\kappa=\min\{\gamma,\frac d2+\frac14\}$, and \eqref{eq:FX1-FX2-d=23}  with $\iota=\kappa=\min\{\gamma,\frac d2+\frac14\}$ to obtain, for $\eta=\min\{\gamma,\frac23\}$
 \begin{align}
 \begin{split}
 &\|A^{\frac\eta2}P
\big(
F(\widetilde{X}_j^{M,N})-F(X_j^{M,N})
\big)
\|
\\
=&
\left\|A^{\frac\eta2}P\int_0^1
F'\big(X(s)+\lambda(\widetilde{X}_j^{M,N})-F(X_j^{M,N})\big)\big(\widetilde{X}_j^{M,N})-F(X_j^{M,N}\big)(1-\lambda)\,\dd \lambda
\right\|
\\
\leq
&C\|e_j^{M,N}\|_1\big(
1
+
\|\widetilde{X}_j^{M,N}\|_\kappa^2
+
\|X_j^{M,N}\|_\kappa^2
\big),
\end{split}
 \end{align}
 where we also used the fact $\eta<\iota=\kappa=\min\{\gamma,\frac d2+\frac14\}$ in dimension one and $\eta<1$ in dimensions two and three.}}
 {\color{black}{Then, by \eqref{eq:sum-x-x^{M,N}-in-H1}, \eqref{eq:boundness-solution-form-full-auxiliary-problem}, \eqref{lem:eq-bound-solution-full-stochatic} and the fact $\widetilde{e}_j^{M,N}\in \dot{H}^1$, one can find that, for $\eta=\min\{\gamma,\frac23\}$  and $\kappa=\min\{\gamma,\frac d2+\frac14\}$,}}
\begin{align}\label{eq:estimate-L2}
\begin{split}
L_2
\leq
&
C
\Big\|
k
\sum_{j=1}^m
t_{m-j+1}^{-\frac{2-\eta}4}
\|A^{\frac\eta2}P
\big(
F(\widetilde{X}_j^{M,N})-F(X_j^{M,N})
\big)
\|
\Big\|_{L^p(\Omega;\mathbb{R})}
\\
\leq
&
C
\Big\|
k
\sum_{j=1}^m
t_{m-j+1}^{-\frac{2-\eta}4}
|\widetilde{e}_j^{M,N}|_1
\,
\big(
1
+
\|\widetilde{X}_j^{M,N}\|_\kappa^2
+
\|X_j^{M,N}\|_\kappa^2
\big)
\Big\|_{L^p(\Omega;\mathbb{R})}
\\
\leq
&
C
\Big\|
\Big(
k
\sum_{j=1}^m
|\widetilde{e}_j^{M,N}|^2_1
\Big)^{\frac 12}
\Big(
k
\sum_{j=1}^m
t_{m-j+1}^{-\frac{2-\eta}2}
\big(
1+\|\widetilde{X}_j^{M,N}\|_\kappa^4
+\|X_j^{M,N}\|_\kappa^4
\big)
\Big)^{\frac 12}
\Big\|_{L^p(\Omega;\mathbb{R})}
\\
\leq
&
C
\Big\|
k
\sum_{j=1}^m
|\widetilde{e}_j^{M,N}|_1^2
\Big\|_{L^p(\Omega;\mathbb{R})}^{\frac12}
\Big\|
k\sum_{j=1}^{m}
t_{m-j+1}^{-\frac{2-\eta}2}
\big(
1+\|\widetilde{X}_j^{M,N}\|_\kappa^4
+\|X_j^{M,N}\|_\kappa^4
\big)
\Big\|_{L^p(\Omega;\mathbb{R})}^{\frac12}
\\
\leq
&
C
(
\lambda_N^{-\frac\gamma2}
+k^{\frac\gamma 4}
).
\end{split}
\end{align}
Finally, inserting  \eqref{eq:A-estimate-full} and \eqref{eq:estimate-L2} into  \eqref{eq:widetilde{e}-decompose-full}
shows \eqref{them:error-estimates-auxiliary-problem} and finishes the proof.
$\square$

\bibliography{bibfile}

\bibliographystyle{abbrv}
 \end{document}